\documentclass[a4paper,10pt]{article}
\usepackage{mathptmx}
\usepackage[T1]{fontenc}
\usepackage[left=1.5cm, right=1.5cm, top=2cm, bottom=2cm, headsep=.5cm]{geometry}
\usepackage{xcolor,amsmath,amssymb,amsthm,mathtools,hyperref,url,subcaption,fancyhdr,caption,appendix,comment}
\hypersetup{colorlinks=true,allcolors=blue}
\usepackage[capitalize,nameinlink,noabbrev]{cleveref}
\theoremstyle{plain}
\newtheorem{theorem}{Theorem}[section]

\newtheorem*{theorem*}{Theorem}

\newtheorem{proposition}[theorem]{Proposition}

\newtheorem{lemma}[theorem]{Lemma}

\newtheorem{corollary}[theorem]{Corollary}

\theoremstyle{remark}
\newtheorem{remark}[theorem]{Remark}

\theoremstyle{definition}

\newcommand{\R}{\mathbf R}
\newcommand{\C}{\mathbf C}
\newcommand{\N}{\mathbf N}
\newcommand{\Q}{\mathbf Q}
\newcommand{\Z}{\mathbf Z}
\newcommand{\T}{\mathbf T}

\newcommand{\zn}{\underline{z}}

\DeclareMathOperator{\diam}{diam}

\title{Limits of Mahler measures in multiple variables} 
\author{François Brunault, Antonin Guilloux, Mahya Mehrabdollahei, Riccardo Pengo}

\date{}

\begin{document}
	
	\maketitle
	
	\begin{abstract}
		We prove that certain sequences of Laurent polynomials, obtained from a fixed multivariate Laurent polynomial $P$ by monomial substitutions, give rise to sequences of Mahler measures which converge to the Mahler measure of $P$.
		This generalises previous work of Boyd and Lawton, who considered univariate monomial substitutions.
		We provide moreover an explicit upper bound for the error term in this convergence, extending work of Dimitrov and Habegger, and a full asymptotic expansion for a family of $2$-variable polynomials, whose Mahler measures were studied independently by the third author. 	
	\end{abstract}
	
	\section{Introduction}
	
	Let $P \in \Z[z]$ be a monic polynomial with integer coefficients. Writing $P(z) = \prod_{j = 1}^d (z - \alpha_j)$ for its complex factorisation, the product $\Delta_k(P) = \prod_{j = 1}^d (\alpha_j^k-1)$ is an integer for every $k \in \N$. The case of the polynomial $P(z)=z-2$ recovers the Mersenne numbers $2^k-1$. In the aim of finding large prime numbers, D.~H.~Lehmer and T.~A.~Pierce \cite{Lehmer_1933, Pierce_1916} developed specific primality tests for these numbers $\Delta_k(P)$. If no root of $P$ lies on the unit circle, then $\lvert \Delta_{k+1}(P)/\Delta_k(P) \rvert$ converges to the real number $M(P) := \prod_{j=1}^d \max(\lvert \alpha_j \rvert,1) \geq 1$. The (logarithmic) Mahler measure of $P$ is then defined as $m(P) = \log M(P)$. In order to maximise the number of prime candidates, Lehmer searched for polynomials with small Mahler measure, and wondered whether there exist irreducible polynomials $P \in \Z[z]$ with $m(P)$ non-zero and arbitrarily small. This seemingly simple question is open to this day, and the Mahler measure of Lehmer's polynomial $z^{10}+z^9-z^7-z^6-z^5-z^4-z^3+z+1$ is still the smallest, non-zero Mahler measure of a polynomial $P \in \Z[z] \setminus \{0\}$ which has been computed (see \cite{Smyth_2008} for a survey).
	
	Thanks to Jensen's formula, $m(P)$ can be expressed as an integral over the unit circle \cite[Proposition 1.4]{Brunault_Zudilin_2020}. This integral generalises naturally to polynomials with several variables, leading to the following general definition of the Mahler measure, which is due to K. Mahler \cite{Mahler_1962}. 
	More precisely, given a $n$-variable Laurent polynomial $P \in \C[z_1^{\pm 1},\dots,z_n^{\pm 1}] \setminus \{0\}$, the Mahler measure of $P$ is defined as
	\[
		m(P) := \int_{[0,1]^n} \log\lvert P(e^{2 \pi i t_1}, \dots, e^{2 \pi i t_n}) \rvert \, d t_1 \cdots d t_n.
	\] 
It can be shown that $m(P) \geq 0$ whenever $P$ has integer coefficients \cite[p.~117]{Boyd_1981}.
	
	One of the most interesting strategies to attack Lehmer's problem has been proposed by Boyd \cite{Boyd_1981_Speculations}. He observed that if the set:
	\[
		\mathcal{M} := \bigcup_{n = 1}^{+\infty} m(\Z[z_1^{\pm 1},\dots,z_n^{\pm 1}] \setminus \{0\}) \subseteq \R_{\geq 0}
	\]
	is closed, then indeed there exists $m_0 > 0$ such that for each $P \in \Z[t]$ either $m(P) = 0$ or $m(P) \geq m_0$. 
	This observation follows from the fact that each Mahler measure $m(P)$ of a multivariate polynomial $P \in \Z[z_1^{\pm 1},\dots,z_n^{\pm 1}] \setminus \{0\}$ is the limit of a sequence of Mahler measures of univariate polynomials. 
	More precisely, Boyd shows that:
	\[
		\lim_{a_1 \to +\infty} \cdots \lim_{a_n \to +\infty} m(P(t^{a_1},\dots,t^{a_n})) = m(P)
	\]
	where each limit is taken independently.
It seems natural to ask what kind of monomial substitutions in the variables of $P$ give the same convergence. Indeed, for a Laurent polynomial $ P \in \C[z_1^{\pm 1},\dots,z_n^{\pm 1}] $ and  a matrix $ A =(a_{i,j})\in \Z^{m \times n} $, one can consider the polynomial $P_A$, in $m$ variables, given by: 
\begin{equation*} 
	P_A(z_1, \ldots, z_m) = P(z_1^{a_{1,1}} \cdots z_m^{a_{m,1}},\dots,z_1^{a_{1,n}} \cdots z_m^{a_{m,n}}).
\end{equation*}	
	The substitutions appearing in the previous limit proven by Boyd are the special case of the row-matrix $A = (a_1,\ldots,a_n)$.
In order to generalise Boyd's result, Lawton \cite{Lawton_1983} considered the quantity $\rho(A)$ associated to the matrix $A$ defined as the smallest $\ell^\infty$-norm of an integer vector in the kernel of $A$:
\[
\rho(A) := \min\{ \lVert v \rVert_\infty \colon v \in \Z^{n} \setminus \{ 0 \}, \ A \cdot v = 0 \}.
\]
Lawton showed that if $A_d$ is a sequence of row-matrices, with $\rho(A_d)\to \infty$, then we have:
\[\lim_{d\to \infty} m(P_{A_d}) = m(P).\]
Recently, Dimitrov and Habegger \cite[Theorem~A.1]{Dimitrov_Habegger_2019} have given an upper bound on the rate of convergence which is a negative power of $\rho(A)$. Strikingly, the exponent depends only on the number of non-vanishing coefficients in $P$. The constant involved depends also on the degree of $P$ and the number $n$ of variables.
	
	Moreover, Smyth \cite{Smyth_2018} used Lawton's result to show that the set $\mathcal{M}$ can be written as a nested ascending union of closed subsets of $\R$.
	In fact, Smyth proves more generally that for every Laurent polynomial $P \in \C[z_1^{\pm 1},\dots,z_n^{\pm 1}] \setminus \{0\}$, the set:
	\begin{equation} \label{eq:Mahler_set}
		\mathcal{M}(P) := \bigcup_{m = 1}^{+\infty} \{ m(P_A) \colon A \in \Z^{m \times n} \ \text{such that} \ P_A \neq 0 \}
	\end{equation}
	is closed.
	Smyth shows moreover that $\mathcal{M}$ is the nested ascending union of the sets $\mathcal{M}\left( \sum_{j = 1}^n z_{2 j - 1} - z_{2 j} \right)$ for $n \to +\infty$ (see \cite[Proposition~14]{Smyth_2018}).
	
	With this context in mind, it seems natural to understand sequences $m(P_{A_d})$ and their convergence when $A_d$ is a sequence of $m\times n$-matrices, and not only of row-matrices.	
	The present paper aims at initialising a systematic study of these sequences. To do so, first of all we devote \cref{sec:higher_Lawton} to the proof of the following theorem (see  \cref{thm:Lawton_generalisation}) which very naturally generalises the theorems of Boyd and Lawton to the multivariate setting:
\begin{theorem} \label{thm:intro}
For every non-zero Laurent polynomial $P \in \C[z_1^{\pm 1}, \ldots, z_n^{\pm 1}] \setminus \{0\}$, and every sequence of integer $m\times n$-matrices $\{A_d\}_{d \in \N} \subseteq \Z^{m \times n}$ such that $\displaystyle\lim_{d \to +\infty} \rho(A_d) = +\infty$, we have that $\displaystyle\lim_{d \to +\infty} m(P_{A_d}) = m(P)$.
\end{theorem}
	
	Further, we obtain in \cref{thm:error-term} an upper bound for the error term $\lvert m(P_A) - m(P) \rvert$, which generalises the bound proved by Dimitrov and Habegger \cite[Theorem~A.1]{Dimitrov_Habegger_2019}.
	In fact, our proof of \cref{thm:error-term}, which occupies the entirety of \cref{sec:error-term}, follows a strategy similar to the one of Dimitrov and Habegger. 
	More precisely, we proceed, as they do, by regularising the function $\log\lvert P \rvert$, and we bound separately the error terms for the regularisations (see \cref{coro:bound2}) and the integrals of the differences between $\log\lvert P \rvert$ and the regularised functions (see \cref{prop:error_regularization}).
	However, our regularisation proceeds by using the smooth functions $\frac{1}{2} \log(\lvert P \rvert^2 + \varepsilon)$, which extend holomorphically to a neighbourhood of the unit torus (see \cref{pro:holCeps}), whereas the regularisation carried out in \cite{Dimitrov_Habegger_2019} uses functions which are not smooth in general.
Let us point out as well that our proof of \cref{prop:error_regularization} relies on an estimate about the volume of the subset of the torus where a polynomial is small (see \cref{thm:improved_DH}), which is a slight improvement on results of Dobrowolski \cite[Theorem~1.3]{Dobrowolski_2017} and Dimitrov-Habegger  \cite[Lemma~A.3]{Dimitrov_Habegger_2019} (see \cref{rmk:Dimitrov_Habegger} for a comparison).
	Along the proof, we show furthermore that if a polynomial $P$ does not vanish on the unit torus $\T^n$, then $m(P_{A})$ tends to $m(P)$ exponentially fast as $\rho(A) \to +\infty$ (see \cref{cor:exponential_convergence_polynomials}).
	
	\cref{sec:speed_of_convergence} gives some insight on which optimal rate of convergence and even on what kind of asymptotic expansion one can expect for the convergence of $m(P_{A})$ towards $m(P)$.
	For the case of $1$-variable Mahler measures converging to $2$-variable ones, the asymptotic expansions of the error term have been studied by Condon \cite{Condon_2012}. His work, which we review in \cref{sec:Condon}, shows that for a large class of $2$-variable polynomials, the rate of convergence for the limit $m(P(z_1,z_1^d)) \to m(P)$ is an integer power of $1/d$, and a full asymptotic expansion for the error term can be obtained. 
	Note that this rate of convergence is much better than the bounds provided by \cite[Theorem~A.1]{Dimitrov_Habegger_2019} and by \cref{thm:error-term}. 
	Moreover, Condon proceeds to give experimental evidences for other polynomials, exhibiting what seems to be a rate of convergence comparable to a rational power of $1/\rho(A)$. The full description of the rate of convergence, even in this particular case, is still open.
	
	However, we exhibit in \cref{sec:Pd} the example of the $4$-variate polynomial $P_\infty(z_1,\dots,z_4) = (1-z_1) (1-z_2) - (1-z_3) (1-z_4)$ and a sequence $A_d$ of $2\times 4$ integer matrices, such that the polynomials $P_{A_d}$ are intimately related to the sequence: 
	\[P_d(z_1,z_2) := \sum_{0 \leq i + j \leq d} z_1^i z_2^j\] and in particular $m(P_{A_d}) = m(P_d)$. 
	The sequence of Mahler measures $m(P_d)$ was thoroughly studied by the third named author of the present paper in \cite{Mehrabdollahei_2021}, where she proved that $m(P_d) \to -18 \cdot \zeta'(-2)$ as $d \to +\infty$. We use \cref{thm:Lawton_generalisation} to give a new proof of this convergence, using the equality $-18 \zeta'(-2) = m(P_\infty)$, which is due to D'Andrea and Lalín \cite[Theorem~7]{D'Andrea_Lalin_2007}. 
	We then provide a complete asymptotic expansion for the error term $m(P_d) - m(P_\infty)$ as $d \to +\infty$ in \cref{thm:Pd_asymptotics}. In particular, we prove that the following asymptotic behaviour holds true when $d \to +\infty$:
	\[
	m(P_d) - m(P_\infty) \sim -\frac{\log(\rho(A_d))}{2 \rho(A_d)^2}.
	\]
The logarithmic term represents a different behaviour than what Condon studied and proved. Moreover, this asymptotic is still much better than our general bound. So, this example shows how far we are from fully understanding the optimal rate of convergence and asymptotic expansion of $m(P_{A_d})$ to $m(P)$ in a general multivariate setting.

\subsection{Historical remarks}

We devote this subsection to a short historical overview of the existing results using and generalising the work of Boyd \cite{Boyd_1981,Boyd_1981_Speculations} and Lawton \cite{Lawton_1983}.
First of all, Boyd himself \cite{Boyd_1981} used an earlier version of this theorem to characterise those Laurent polynomials $P \in \Z[\zn^{\pm 1}] \setminus \{0\}$ such that $m(P) = 0$. 
Moreover, Lawton's result has been used by Schinzel \cite{Schinzel_1997} to provide an explicit bound on the Mahler measure of a polynomial, which generalises a classical result of Gonçalves \cite{Goncalves_1950} (see also \cite[Theorem~1.22]{Everest_Ward_1999}).  
Furthermore, the work of Boyd and Mossinghoff \cite{Boyd_Mossinghoff_2005}, later generalised by Otmani, Rhin and Sac-Épée \cite{Otmani_Rhin_Sac-Epee_2019}, used Lawton's result as a starting point for an investigation of the genuine limit points in the set $\mathcal{M}$. On the other hand, Dobrowolski \cite{Dobrowolski_2012} used Lawton's limit formula to answer a question of Schinzel.
Moving on, Dubickas and Jankauskas \cite{Dubickas_Jankauskas_2013} used Lawton's theorem to construct many non-reciprocal univariate polynomials whose Mahler measures lie in the interval $[m(z_1^3 - z_1 - 1),m(1+z_1+\dots+z_n)]$, whereas Dobrowolski and Smyth \cite{Dobrowolski_Smyth_2017}, as well as Akhtari and Vaaler \cite{Akhtari_Vaaler_2019}, used the theorem of Lawton to study Mahler measures of polynomials with a bounded number of monomials.
Finally, Dubickas \cite{Dubickas_2018} and Habegger \cite{Habegger_2018} used Lawton's result in their investigations of sums of roots of unity, whereas, as we already mentioned, Smyth \cite{Smyth_2018} used Lawton's limit formula to prove that the sets $\mathcal{M}(P)$ defined in \eqref{eq:Mahler_set} are closed.

Let us point out that Lawton's result has found applications also outside number theory.
First of all, Lind, Schmidt and Ward \cite{Lind_Schmidt_Ward_1990} used it to provide a lower bound for the entropy of the dynamical system associated to a Laurent polynomial in terms of its Mahler measure.
Moreover, the work of Silver and Williams \cite{Silver_Williams_2004,Silver_Williams_2012}, later generalised by Raimbault \cite{Raimbault_2012} and Lê \cite{Le_2014}, applied Lawton's result to knot theory, in order to study the convergence of Mahler measures of Alexander polynomials, and the growth of homology under surgery operations. Staying in the realm of knot theory, the work of Champanerkar and Kofman \cite{Champanerkar_Kofman_2005,Champanerkar_Kofman_2006}, later generalised by Cai and Todd \cite{Cai_Todd_2014}, used Lawton's theorem to study Mahler measures of Jones polynomials.
Moving to the world of von Neumann algebras, Deninger \cite[Theorem~17]{Deninger_2009} proved a continuity result for Fuglede-Kadison determinants on the space of marked groups, which implies Lawton's result under the strong assumption that $P$ does not vanish on the torus $\T^n$.
Note that Deninger's result is reminiscent of the classical theorem of Szegö \cite[Theorem~1]{Deninger_2009}, which approximates univariate Mahler measures in terms of Toeplitz determinants. A multivariate analogue of Szegö's result has been recently found by Hajli \cite{Hajli_2020}.
Finally, Lawton's result has been used by Lück in functional analysis, to study spectral density functions \cite{Luck_2015} and twists of $L^2$-invariants \cite{Luck_2018}.

To conclude this subsection, let us mention some existing generalisations and improved versions of Lawton's result. First of all, Champanerkar and Kofman prove in \cite[Lemma~3.3]{Champanerkar_Kofman_2006} that one can perform signed monomial substitutions. 
Moreover, Duke \cite[Theorem~6]{Duke_2007} provides the first term in the asymptotic expansion of the difference $m(z_1^n + z_1^m + 1) - m(z_1 + z_2 + 1)$ as $(n,m) \to +\infty$. 
Furthermore, Lalín and Sinha \cite{Lalin_Sinha_2011} mention a generalisation of Lawton's theorem \cite[Theorem~30]{Lalin_Sinha_2011} to the multiple Mahler measure, introduced in previous work of Kurokawa, Lalín and Ochiai \cite{Kurokawa_Lalin_Ochiai_2008}. 
Such a generalisation was rigorously proved by Issa and Lalín \cite{Issa_Lalin_2013}, who dealt also with the generalised Mahler measures defined by Gon and Oyanagi \cite{Gon_Oyanagi_2004}. 
On the other hand, Carter, Lalín, Manes, Miller and Mocz \cite[Proposition~1.3]{Carter_Lalin_Manes_Miller_Mocz_2021} recently proved a weak generalisation of Lawton's result to dynamical Mahler measures (introduced in \cite[Definition~1.1]{Carter_Lalin_Manes_Miller_Mocz_2021}).

In addition to the previously mentioned results, Dobrowolski \cite{Dobrowolski_2017} generalised a crucial estimate of Lawton \cite[Theorem~1]{Lawton_1983} on the measure of the set of points $\zn \in \T^n$ where a polynomial is small. Similar bounds have been provided by Lück \cite[Proposition~2.1]{Luck_2015}, Habegger \cite[Lemma~A.4]{Habegger_2018}, and Dimitrov and Habegger \cite[Lemma~A.3]{Dimitrov_Habegger_2019} (see \cref{rmk:Luck,rmk:Dimitrov_Habegger} for a comparison).
Finally, Gu and Lalín \cite[Proposition~8]{Gu_Lalin_2021} have recently proved a multivariate convergence of Mahler measures, which can be obtained as a corollary of \cref{thm:intro}, for one particular family of polynomials (see for more details \cref{rmk:Gu_Lalin}).

\section{Notation and conventions}

For the reader's convenience, we collect in this section the notation that we most frequently use in the rest of the paper.

\subsection{Generalities}

We let $\N = \{0,1,\dots\}$ denote the natural numbers, $\Z$ denote the integers, $\R$ denote the real numbers, $\C$ denote the complex ones, and $\C^\times = \C\setminus\{0\}$.
For any $n \in \N$, we denote by $\underline{z}_n = (z_1,\dots,z_n)$ the coordinates of $\C^n$, and by 
\[
	\T^n := \{ \underline{z}_n \in \C^n \colon \lvert z_1 \rvert = \cdots = \lvert z_n \rvert = 1 \}
\] 
the $n$-dimensional real-analytic unit torus.
Moreover, for every $p \in \R_{\geq 1}$, we let $\lVert \cdot \rVert_p \colon \C^n \to \R_{\geq 0}$ denote the $\ell^p$-norm, defined for every $\zn_n \in \C^n$ by: 
\[
\lVert \zn_n \rVert_p := \left( \sum_{j = 1}^{n} \lvert z_j \rvert^p \right)^{1/p}
\]
and we let $\lVert \cdot \rVert_\infty \colon \C^n \to \R_{\geq 0}$ denote the $\ell^\infty$-norm, defined by $\lVert \zn_n \rVert_\infty := \max\left\{ \lvert z_1 \rvert,\ldots,\lvert z_n \rvert \right\}$.
Finally, for any natural number $n$ and any real number $\delta > 0$, we define the \emph{annulus} 
\begin{equation} \label{not:C_delta}
	\mathcal{C}_\delta := \left\{ \zn_n \in (\C^\times)^n \colon \sum_{i = 1}^n \lvert \log\lvert z_i \rvert \rvert \leq \delta \right\}
\end{equation}
which is a closed neighbourhood of the torus $\T^n$ in $(\C^\times)^n$.

\subsection{Matrices}
\label{not:matrices}
Fix a matrix $A = (a_{i,j}) \in \Z^{m \times n} $. We denote by $d(A)$ the dimension of the linear subspace $\ker(A) \subseteq \R^n$, and by $\Lambda_A := \ker(A) \cap \Z^n$ the integer lattice within this subspace. Moreover, we introduce the quantity:
\begin{equation*} 
	\rho(A) := \min\{ \lVert v \rVert_\infty \colon v \in \Lambda_A, \, v \neq 0 \}
\end{equation*}
which is the first successive minimum of the lattice $\Lambda_A$ with respect to the $\ell^\infty$-norm. By convention, we put $\rho(A) = +\infty$ when $\Lambda_A = \{0\}$.
Finally, we consider the monomial substitution:
\begin{equation*}
	\zn_m^A := (z_1^{a_{1,1}} \cdots z_m^{a_{m,1}},\dots,z_1^{a_{1,n}} \cdots z_m^{a_{m,n}})
\end{equation*}
and for any Laurent polynomial $P \in \C[z_1^{\pm 1},\dots,z_n^{\pm 1}]$, we define $P_A \in \C[z_1^{\pm 1},\dots,z_m^{\pm 1}]$ by setting $P_A(z_1,\dots,z_m) := P(\zn_m^A)$. In particular, we consider vectors $v = (v_1,\dots,v_n) \in \Z^n$ as column matrices, so that $\zn_n^v = z_1^{v_1} \cdots z_n^{v_n}$ is a monomial.

\subsection{Measure theory}
For every $n \in \N$, we denote by $\mu_n := \frac{1}{(2 \pi i)^n} \frac{d z_1}{z_1} \wedge \cdots \wedge \frac{d z_n}{z_n}$ the probability Haar measure on $\T^n$.
More generally, for every matrix $A = (a_{i,j}) \in \Z^{m \times n}$, we let $\mu_A$ be the probability measure on $\T^n$ defined as the push-forward of $\mu_m$ along the map $\T^m \to \T^n$ given by $\zn_m \mapsto \zn_m^A$. 
Note in particular that $\mu_{\operatorname{Id}_n} = \mu_n$. Finally, for every non-zero Laurent polynomial $P \in \C[\zn_n^{\pm 1}] \setminus \{0\}$, we let:
\begin{equation*}
	m(P) := \int_{\T^n} \log\lvert P(\zn_n) \rvert d\mu_n(\zn_n) \in \R
\end{equation*}
denote the \textit{logarithmic Mahler measure} of $P$.

\subsection{Fourier coefficients}
\label{not:Fourier_coefficients}
For every integrable function $f \colon \T^n \to \C$, and every vector $v \in \Z^n$, we denote by:
\[
c_v(f) := \int_{\T^n} \frac{f(\zn_n)}{\zn_n^v} d\mu_n(\zn_n)
\]
the corresponding Fourier coefficient. In particular, if $P \in \C[z_1^{\pm 1},\dots,z_n^{\pm 1}]$ is a Laurent polynomial, then $P(\zn_n) = \sum_{v \in \Z^n} c_v(P) \cdot \zn_n^v$.

\subsection{Polynomials}
\label{not:polynomials}
Fix a non-zero Laurent polynomial $P \in \C[z_1^{\pm 1},\dots,z_n^{\pm 1}]$. We denote by $N_P \subseteq \R^n$ its \textit{Newton polytope}, which is the convex hull in $\R^n$ of the \textit{support} set $\operatorname{supp}(P) := \{ v \in \Z^n \colon c_v(P) \neq 0 \}$. We also denote by $\diam(P)$ the \textit{diameter} of $P$, which is the smallest $d \in \N$ such that $N_P$ is contained inside a translate of $[0,d]^n$.

We define $k(P)$ to be the number of non-zero coefficients $P$, and for every $1 \leq i \leq n$, we define $k_i(P)$ to be the number of non-zero coefficients of $P$ seen as a polynomial in $z_i$. We also denote $\kappa(P) = \max_{1 \leq i \leq n}(k_i(P))$.

We write $L_1(P) := \sum_{v \in \Z^n} \lvert c_v(P) \rvert$ for the \textit{length} of $P$, and $L_\infty(P) := \max_{v \in \Z^n} \lvert c_v(P) \rvert$ for the \textit{modulus} of $P$.

Furthermore, we let $V_P \hookrightarrow \mathbf{G}_m^n$ be the hypersurface defined by $P$, so that $V_P(\C) := \{ \zn_n \in (\C^\times)^n \colon P(\zn_n) = 0 \}$. We also define the \textit{conjugate reciprocal} of $P$ by:
\begin{equation*}
	P^\ast(\zn_n) := \overline{P}(\zn_n^{- 1}) = \sum_{v \in \Z^{n}} \overline{c_v(P)} \cdot \zn_n^{-v}.
\end{equation*}

Finally, for every $r \geq 0$ we define the set $S(P,r) := \{ \zn \in \T^n \colon \lvert P(\zn) \rvert \leq r \} \subseteq \T^n$.

\subsection{Constants}

For every non-zero Laurent polynomial $P \in \C[\zn_n^{\pm 1}] \setminus \{0\}$ we define a constant: 
\begin{equation} \label{not:rho_0}
\rho_0(P) := \max\left\{\diam(P)+1, \; 7  \diam(P)^2, \; \exp(2 (k-1) \max(n,5)) \right\} \in \R_{> 0},
\end{equation}
and a further family of constants: 
\begin{equation} \label{not:a_eps}
	\delta_\varepsilon(P) := \min\left( \frac{\sqrt{\varepsilon}}{\diam(P) L_1(P)}, \frac{\log(4/3)}{\diam(P)} \right) \in \R_{> 0}
\end{equation}
depending on a positive real number $\varepsilon > 0$. 

\section{A higher dimensional analogue of Lawton's theorem}
	\label{sec:higher_Lawton}
	
	The aim of this section is to show that the Mahler measure $m(P) := \int_{\T^n} \log\lvert P(\zn_n) \rvert d\mu_n(\zn_n)$ of any non-zero Laurent polynomial $P(\zn_n) \in \C[\zn_n^{\pm 1}] \setminus \{0\}$ can be approximated by suitable sequences of ``lower-dimensional'' Mahler measures, as specified in the following theorem (see also \cref{thm:intro}).
		
	\begin{theorem} \label{thm:Lawton_generalisation}
		Let $ n \in \N $ be an integer, and $P(\zn_n) \in \C[\zn_n^{\pm 1}] \setminus \{0\}$ be a non-zero Laurent polynomial. Then, for every sequence of matrices $A_d \in \Z^{m_d \times n}$ such that $\displaystyle\lim_{d \to +\infty} \rho(A_d) = +\infty$, we have the convergence $\displaystyle\lim_{d \to +\infty} m(P_{A_d}) = m(P)$.
	\end{theorem}
	
	\subsection{Convergence of measures and integrals}
	
	In order to prove \cref{thm:Lawton_generalisation}, we start by relating the growth of $\rho(A)$ to the weak convergence of the push-forward measures $\mu_A$:
	
	\begin{lemma} \label{lem:matrix_measure_convergence}
		Fix $n \in \Z_{\geq 1}$, and let $A_d \in \Z^{m_d \times n}$ be a sequence of integral matrices, with fixed number of columns, such that $\rho(A_d) \to +\infty$ as $d \to +\infty$. Then the sequence of measures $\mu_{A_d}$ on $\T^n$ converges weakly to the measure $\mu_{\operatorname{Id}_n}$.
	\end{lemma}
	\begin{proof} 
		This result is classical. We follow the lines of \cite[Lemma 1]{Boyd_1981}, which treats the case when $m_d = 1$ for every $d$. By the definition of weak convergence and push-forward of measures, and by Weierstra{\ss} approximation, it is sufficient to prove that:
		\begin{equation} \label{eq:weak_convergence_polynomials}
			\lim_{d \to +\infty} \left( \int_{\T^{m_d}} Q(\zn_{m_d}^{A_d}) \, d \mu_{m_d}(\zn_{m_d}) \right) = \int_{\T^n} Q(\zn_n) \, d \mu_{n}(\zn_n)
		\end{equation}
		for every Laurent polynomial $Q(\zn_n) \in \C[\zn_n^{\pm 1}]$.
		We see now immediately that for every $d \in \N$, the following identities hold true:
		\begin{equation} \label{eq:Fourier_expansion_polynomials}
			\int_{\T^{m_d}} Q(\zn_{m_d}^{A_d}) \, d \mu_{m_d}(\zn_{m_d}) = \sum_{v \in \Z^{n}} c_v(Q) \cdot \int_{\T^{m_d}} \zn_{m_d}^{A_d \cdot v} \ d \mu_{m_d}(\zn_{m_d}) = \sum_{\substack{v \in \Z^{n} \\ A_d \cdot v = 0}} c_v(Q) = \sum_{v\in \Lambda_{A_d}} c_v(Q)
		\end{equation}
		Now, set $R := \max\{\lVert v \rVert_\infty \colon v \in \mathrm{supp}(Q) \}$. If $\rho(A_d) > R$, then the only vector $v \in \Lambda_{A_j}$ for which it may happen that $c_v(Q)\neq 0$ is the null vector $v = 0$. In this case, we have the identity:
		\[
		\sum_{v\in \Lambda_{A_d}} c_v(Q) = c_0(Q) = \int_{\T^n} Q(\zn_n) \, d \mu_{\operatorname{Id}_n}(\zn_n)
		\]
		which, combined with \eqref{eq:Fourier_expansion_polynomials}, shows \eqref{eq:weak_convergence_polynomials}, because the sequence on the left is eventually constantly equal to the right hand side.
	\end{proof}
	
Weak convergence of measures implies the convergence of integrals of any bounded continuous function. Unfortunately, we would like a convergence of integrals of $\log|P|$, which is singular. However, uniform estimates on $L^2$-norms are enough to guarantee that the weak-convergence of measures implies convergence of integrals, as shown in the following general \cref{lem:convergence of integrals}. In this lemma, we choose to work with continuous functions possibly having $+\infty$-values, for which the integral for any measure on the torus is naturally defined (possibly $+\infty$), as explained for instance in \cite[Chapter 1]{Rudin}. 
	\begin{lemma}\label{lem:convergence of integrals}
		Let $\nu_k$ be a sequence of probability measures on $\T^n$, which converges weakly to some probability measure $\nu_\infty$. 
		Let $f \colon \T^n \to \R\cup\{+\infty\}$ be a continuous function, which is uniformly $L^2$ for the family $\{ \nu_k \colon k \in \N \cup \{\infty\} \}$.
		Then we have the convergence $\int_{\T^n} fd\nu_k \to \int_{\T^n}f d\nu_\infty$ as $k \to +\infty$.
	\end{lemma}
	\begin{proof}
		By assumption, there exists a positive real number $C \in \R_{> 0}$ such that $\int_{\T^n} \lvert f \rvert^2 d \nu_k \leq C$ for every $k \in \N \cup \{\infty\}$.
		Fix $\varepsilon >0$ and let $\lambda = \frac{C}{\varepsilon}$. Define the set $S_\lambda = \{t\in \T^n \colon \lvert f(t) \rvert > \lambda\}$. The $L^2$-bounds yield, for any $k \in \N \cup \{\infty\}$:
		\begin{equation}\label{eq:L2}
			\left|\int_{S_\lambda} (f-\lambda) d\nu_k \right| \leq 2\int_{S_\lambda} \left|f \right|d\nu_k \leq 2 \int_{S_\lambda} |f|\frac{|f|}{\lambda}d\nu_k \leq \frac{2}{\lambda} \int_{\T^n} |f|^2d\nu_k\leq \frac{2C}{\lambda}=2\varepsilon.
		\end{equation}
		Now, let $\tilde f$ be the continuous function $\min(f,\lambda)$, which is bounded from above by $\lambda$. For every $\zn \in \T^n$, we have the equality $f(\zn) = \tilde{f}(\zn) + (f(\zn)-\lambda) \cdot \chi_{S_\lambda}(\zn)$, where $\chi_{S_\lambda}$ denotes the characteristic function of $S_\lambda$.
		Hence, for all $k \in \N$ we have the bound:
		\begin{equation} \label{eq:L2_bound}
			\left|\int_{\T^n} fd\nu_k -\int_{\T^n} fd\nu_\infty \right|  \leq \left| \int_{\T^n} {\tilde f} d\nu_k - \int_{\T^n} {\tilde f} d\nu_\infty \right| + \left| \int_{S_\lambda} (f-\lambda)d\nu_k\right| + \left| \int_{S_\lambda} (f-\lambda)d\nu_\infty\right|.
		\end{equation}
		The last two terms on the right hand side of \eqref{eq:L2_bound} are bounded by $2\varepsilon$, thanks to \eqref{eq:L2}. For $k$ big enough, by the convergence $\nu_k\to \nu_{\infty}$, the first one is less than $\varepsilon$. So we have proven that, for $k$ big enough, we have:
		\[
			\left|\int_{\T^n} f d\nu_k -\int_{\T^n} f d\nu_\infty \right|  \leq 5\varepsilon
		\]
		which shows that $\int_{\T^n} fd\nu_k \to \int_{\T^n}f d\nu_\infty$ as $k \to +\infty$.
	\end{proof}	
	
\subsection{Uniform \texorpdfstring{$L^2$}{L2}-bounds and convergence of Mahler measures}

Our goal is to prove that $m(P_A) = \int_{\T^n}\log|P| d\mu_A$ converges to $m(P) = \int_{\T^n}\log|P| d\mu_A$. From the previous results we know that an uniform $L^2$-bound for this functions would grant the convergence. The following estimate is essentially obtained by Dimitrov and Habegger in \cite[Appendix A]{Dimitrov_Habegger_2019}, where they deal with Lawton theorem and improves the rate of convergence (see also \cite{Habegger_2018}).

\begin{proposition}[Dimitrov \& Habegger]\label{cor:uniform L2}
		Let $ n, k \in \N $ be two integers. Then, there exists a constant $C > 0$ such that, for every non-zero Laurent polynomial $P(\zn_n) \in \C[\zn_n^{\pm 1}] \setminus \{0\}$ with $k(P) = k$ and $L_\infty(P) = 1$, and every matrix $A \in \Z^{m \times n}$ with $\rho(A)>\diam(P)$ and $m \leq n$, the following holds:
		$$\|\log|P|\|^2_{2,\mu_A} := \int_{\T^n} \left| \log|P(\zn_m^A)|\right|^2 d\mu_A \leq C \quad{ and } \quad \|\log|P|\|^2_{2,\mu_n}\leq C.$$
\end{proposition}
\begin{proof}
	
	A direct computation, similar to the one carried out in \eqref{eq:Fourier_expansion_polynomials}, proves that, for any matrix $A \in \Z^{m \times n}$, we have:
	\begin{equation*}
		P(\zn_m^A) = \sum_{v \in \Z^{n}} c_v(P) \zn_m^{A \cdot v} = \sum_{w \in \Z^m} \left( \sum_{\substack{v \in \Z^n \\ A \cdot v = w}} c_v(P) \right) \zn_m^w.
	\end{equation*}
	Two vectors $v, v' \in \Z^n$ contribute non-trivially to the same monomial in the above sum if and only if $c_v(P)\neq 0$,  $c_{v'}(P)\neq 0$ and $A \cdot v = A \cdot v'$, or equivalently $v-v' \in \Lambda_A$. By definition of $\diam(P)$ (see \cref{not:polynomials}), in this case, we have $\lVert v-v' \rVert_\infty \leq \diam(P)$. 
	We see that if $\rho(A)>\diam(P)$ the only possibility is $v - v' = 0$. In other terms, each monomial of $P_A$ comes from a single monomial of $P(\zn_m^A)$ with the same coefficient and no compensations.
		
	So, for any $A$ with $\rho(A)>\diam(P)$, the polynomial $P_A = P(\zn_m^A)$ has $m\leq n$ variables, $k(P_A) = k(P) = k$ non-vanishing coefficients, and $L_\infty(P_A) = L_\infty(P) = 1$. Our proposition comes then directly from the estimates of Dimitrov and Habegger. They show in \cite[Lemma A.3(i)]{Dimitrov_Habegger_2019}
	that for every $l,k \in \Z_{\geq 1}$, there exists a constant $C_{l,k}>0$ such that for any Laurent polynomial $Q \in \C[\zn_l^{\pm 1}]$, with $k(Q) = k$ and $L_\infty(Q) = 1$, we have:
	\[
	\int_{\T^l} (\log|Q|)^2 d\mu_l \leq C_{l,k}.
	\]
	Thanks to the considerations in the previous paragraph, this bound applies both to $Q = P$ and to $Q = P_A$, hence we can take $C := \max\{ C_{m,k} \colon m \leq n \}$.
	\end{proof}
	
\cref{thm:Lawton_generalisation} is now an easy consequence of the other results proved in this section:
	
	\begin{proof}[Proof of \cref{thm:Lawton_generalisation}]
		
		We first make an easy reduction: up to multiplying $P$ by a constant $a$, we may and will assume that $L_\infty(P) = 1$. Indeed, we have, for all $a\in \C^*$, both $m(aP)=\log|a|+m(P)$ and $m(aP_{A_j}) = \log|a| + m(P_{A_j})$. So the problem of convergence is equivalently solved for $P$ or $aP$. Observe moreover that, for every $d \in \N$, we have the following identities: 
		\[
			m(P_{A_d}) = \int_{\T^{m_d}} \log|P(\zn_{m_d}^{A_d})|d\mu_{m_d} = \int_{\T^n}\log|P|d\mu_{A_d}.
		\]
		
		Let $d_0 \in \N$ be any natural number such that $\rho(A_d)\geq \diam(P)$ for every $d \geq d_0$. From \cref{cor:uniform L2}, we know that the function $\log|P|$ is uniformly $L^2$ for the family $\{\mu_{A_d}, d\geq d_0\} \cup \{\mu_n\}$.
		Moreover, we know from \cref{lem:matrix_measure_convergence} that the family $\mu_{A_d}$ converges weakly to $\mu_n$ as $d \to +\infty$. Thus, we have:
		\[
			\lim_{d \to +\infty} m(P(\zn_{m_d}^{A_d})) = \lim_{d \to +\infty} \int_{\T^n}\log|P|d\mu_{A_d} = \int_{\T^n}\log|P|d\mu_n = m(P)
		\]
		thanks to \cref{lem:convergence of integrals}.
	\end{proof}

We will see an example of application of \cref{thm:Lawton_generalisation} in \cref{sec:Pd}. 
Meanwhile, we will devote the following section to a more careful analysis of the convergence, which will provide an upper bound on its rate.

\section{An error term in the convergence}
\label{sec:error-term}

The aim of this section is to improve \cref{thm:Lawton_generalisation} by providing an explicit upper bound for the error term $\lvert m(P_A) - m(P) \rvert$, where $P \in \C[\zn_n^{\pm 1}] \setminus \{0\}$ is a non-zero Laurent polynomial, and $A \in \Z^{m \times n}$ is an integral matrix. We will assume without loss of generality that $P$ is not a monomial (\textit{i.e.} $k(P) > 1$), and that $n \geq 2$, because otherwise $m(P_A) = m(P)$ for every non-zero integral matrix $A \in \Z^{m \times n}$. Then, we obtain the following result, which generalises \cite[Theorem~A.1]{Dimitrov_Habegger_2019} to higher dimensions:

\begin{theorem} \label{thm:error-term}
	Fix two integers $k,n \geq 2$. Let $P \in \C[\zn_n^{\pm 1}]$ be a Laurent polynomial with $k(P) = k$ non-zero coefficients, and let $\rho_0(P)$ be the constant defined in \eqref{not:rho_0}. Then, for every $m \in \Z_{\geq 1}$ and every matrix $A \in \Z^{m \times n}$ such that $\rho(A) \geq \rho_0(P)$, the following inequality holds:	
	\begin{equation*} 
	\lvert m(P_A)-m(P) \rvert \leq 8 \cdot (36ek)^{n-1} \cdot \log(\rho(A))^{n} \left( \frac{\diam(P)}{\rho(A)} \right)^\frac{1}{k-1}.
	\end{equation*}
\end{theorem}

\subsection{An explicit exponential convergence for polynomials without toric points}

The aim of this section is to show that, for a Laurent polynomial $P \in \C[\zn_n^{\pm 1}]$ which does not vanish on the unit torus $\T^n$, the convergence $m(P_A) \to m(P)$ as $\rho(A) \to +\infty$ is exponentially fast, and its speed can be explicitly bounded, as we will see in \cref{cor:exponential_convergence_polynomials}. 
This result follows easily from the more general \cref{thm:exponential_convergence_functions}, which we will use later on in the proof of \cref{coro:bound2}, which in turn plays a crucial part in the proof of \cref{thm:error-term}. 
The proof of \cref{thm:exponential_convergence_functions} uses crucially the standard fact that the Fourier coefficients of a holomorphic function on a neighbourhood of the torus $\T^n$ decay exponentially.

\begin{theorem} \label{thm:exponential_convergence_functions}
	Fix two natural numbers $n,m \geq 1$, an open $U \subseteq \C^n$ containing $\T^n$, and a holomorphic function $f \colon U \to \C$. 
	Then, for every real number $\delta > 0$ such that $U$ contains the annulus $\mathcal{C}_\delta$ defined in \eqref{not:C_delta}, and every matrix $A \in \Z^{m \times n}$ such that $\rho(A) \geq  \frac{2 d(A)}{3 \delta}$, the following estimate holds:
	\begin{equation*}
	\left\lvert \int_{\T^n} f(\zn) d \mu_A(\zn) - \int_{\T^n} f(\zn) d \mu_n(\zn) \right\rvert \leq (d(A) + 1) 3^{d(A)} \cdot \frac{\max_{\mathcal{C}_\delta} \lvert f \rvert}{\exp(\delta \rho(A))}
	\end{equation*}
	where $d(A) := \dim(\ker(A))$.
\end{theorem}
\begin{proof}
	Let $r=e^\delta$ and $\Lambda_A = \ker(A) \cap \Z^n$, as in \cref{not:matrices}.
	For any $v \in \Z^n$, write $c_v(f) := \int_{\T^n} f(\zn) \zn^{-v} d\mu_n(\zn)$ for the $v$-th Fourier coefficient of $f$, as in \cref{not:Fourier_coefficients}.
	Since $f$ is holomorphic on $U$, the Fourier series $\sum_{v \in \Z^n} c_v(f) \zn^v$ converges normally to $f$ on $\T^n$, and the dominated convergence theorem gives
	\begin{equation} \label{eq:integral_difference}
		\begin{aligned}
		\int_{\T^n} f(\zn) d \mu_A(\zn) - \int_{\T^n} f(\zn) d \mu_n(\zn) & = \left( \sum_{v \in \Z^n} c_v(f) \int_{\T^n} \zn^{v} d\mu_A(\zn) \right) - c_0(f) = \left( \sum_{v \in \Lambda_A} c_v(f) \right) - c_0(f) = \sum_{v \in \Lambda_A \backslash \{0\}} c_v(f).
		\end{aligned}
	\end{equation}
	To bound the Fourier coefficients $c_v(f)$, we use the holomorphicity of $f$ on $U$. To be more precise, let us associate to every vector $\underline{h} = (h_1,\dots,h_n) \in \R^n$ the torus $\T_{\underline{h}} := \{ \zn \in \C^n \colon \lvert z_j \rvert = e^{h_j}, \ \forall j \in \{1,\dots,n\} \}$. Then, for every $v \in \Z^n \backslash \{0\}$ and every $\underline{h} \in \R^n$ such that $\lVert \underline{h} \rVert_1 \leq \delta$, the homotopy invariance of integrals of holomorphic functions implies
	\begin{equation} \label{eq:Fourier_homotopy}
		c_v(f) = \int_{\T_{\underline{h}}} f(\zn) \zn^{-v} d\mu_n(\zn)
	\end{equation}
	because $\T_{\underline{h}} \subseteq \mathcal{C}_\delta \subseteq U$ by assumption.
	Now let $j_0 \in \{1,\dots,n\}$ be any integer such that $\lVert v \rVert_\infty = \lvert v_{j_0} \rvert$, and take $\underline{h} \in \R^n$ to be the vector with $h_j := 0$ for any $j \in \{1,\dots,n\} \setminus \{j_0\}$, and $h_{j_0} := \delta \cdot v_{j_0}/\lvert v_{j_0} \rvert$. Then, using \eqref{eq:Fourier_homotopy} we see that:
	\begin{equation} \label{eq:bound_fourier_coefficient}
	|c_v(f)| \leq \max_{\mathcal{C}_\delta} \lvert f \rvert \cdot r^{- \lVert v \rVert_\infty}.
	\end{equation}
	Combining \eqref{eq:integral_difference} and \eqref{eq:bound_fourier_coefficient}, we get
	\begin{equation} \label{eq:bound_holomorphic_integral_difference}
		\left\lvert \int_{\T^n} f(\zn) d \mu_A(\zn) - \int_{\T^n} f(\zn) d \mu_n(\zn) \right\rvert \leq \sum_{v \in \Lambda_A \backslash \{0\}} \lvert c_v(f) \rvert \leq \max_{\mathcal{C}_\delta} \lvert f \rvert \sum_{v \in \Lambda_A \backslash \{0\}} r^{- \lVert v \rVert_\infty}.
	\end{equation}
The only remaining step to prove the theorem is to bound the sum appearing in the right-hand-side of \eqref{eq:bound_holomorphic_integral_difference}. It is an independent estimate, which we state separately in \cref{lem:estimate_lattice_series}. Note that we fulfil its assumptions: $\Lambda_A$ is a lattice of full rank inside the vector space $\ker(A)$ of dimension $d(A)$. Moreover, its first successive minimum with respect to the $\ell^\infty$-norm is by definition $\rho(A)$. Eventually, we have by assumption $\rho(A)\log(r) = \rho(A)\delta \geq 2d(A)/3$, so we can use the bound \eqref{eq: bound sum lattice 2} of the following \cref{lem:estimate_lattice_series}.
\end{proof}

\begin{lemma} \label{lem:estimate_lattice_series}
	Fix a real vector space $V$ of finite dimension $d \in \Z_{\geq 1}$, and a norm $\lVert \cdot \rVert \colon V \to \R_{\geq 0}$. Let $\Lambda \subseteq V$ be a lattice of full rank, and denote by 
$\rho$ the norm of the shortest vector of $\Lambda \setminus \{0\}$ with respect to $\lVert \cdot \rVert$.
	Then, we have the following estimate:
	\begin{equation} \label{eq: bound sum lattice 1}
		\sum_{v \in \Lambda\setminus \{0\}} r^{- \lVert v \rVert} \leq \frac{3^d d!}{r^\rho} \sum_{k=0}^d \frac{1}{(d-k)!} \Bigl(\frac{2}{3\rho \log(r)}\Bigr)^k \qquad (r>1).
	\end{equation}
In particular, if $\rho \log(r) \geq 2d/3$, we have
	\begin{equation} \label{eq: bound sum lattice 2}
		\sum_{v \in \Lambda\setminus \{0\}} r^{- \lVert v \rVert} \leq \frac{(d+1) 3^d}{r^\rho} \qquad (r>1).
	\end{equation}
\end{lemma}
\begin{proof}
	First of all, set $B(x,q) := \{ y \in V \colon \lVert y - x \rVert \leq q \}$ and $N_q := \lvert B(0,q) \cap \Lambda \rvert$ for every $x \in V$ and $q \geq 0$. Now, observe that: 
	\begin{equation} \label{eq:lattice_points_first_inequality}
		\begin{aligned}
		\sum_{v \in \Lambda\setminus \{0\}} r^{- \lVert v \rVert} &= \sum_{q = 1}^{+\infty} \lvert \{ v \in \Lambda \colon \lVert v \rVert = q \} \rvert \cdot r^{-q} = \sum_{q = 1}^{+\infty} (N_q - N_{q - 1}) \cdot r^{-q} = \\
		&= -\frac{N_{\rho - 1}}{r^{\rho - 1}} + \log(r) \int_{\rho - 1}^{+\infty} \frac{N_{\lfloor t \rfloor}}{r^t} dt = -\frac{1}{r^\rho} + \log(r) \int_{\rho}^{+\infty} \frac{N_{\lfloor t \rfloor}}{r^t} dt
		\end{aligned}
	\end{equation}
	as follows from Abel's summation formula (see \cite[Theorem~4.2]{Apostol_1976}), since $N_q = 1$ if $q \leq \rho - 1$.
	Moreover, note that:
	\[
		 \bigsqcup_{x \in B(0,q) \cap \Lambda} B\left(x,\frac{\rho}{2} \right)^\circ \subseteq B(0,q + \rho/2) \qquad (q \geq 0)
	\]
	where $B(x,\rho/2)^\circ := \{y \in V \colon \lVert y - x \rVert < \rho/2 \}$.
	This inclusion, together with the fact that $\mathrm{vol}(B(x,\rho/2)) = \mathrm{vol}(B(x,\rho/2)^\circ)$, provides the bound
	\begin{equation} \label{eq:bound_integral_points}
		N_q \leq \frac{\mathrm{vol}\left(B\left(0,q + \frac{\rho}{2}\right)\right)}{\mathrm{vol}\left(B\left(0,\frac{\rho}{2}\right)\right)} = \left( \frac{2 q}{\rho} + 1 \right)^d \qquad (q \geq 0)
	\end{equation}
	 which is also proved in \cite[Theorem~2.1]{Betke_Henk_Wills_1993}. Applying the bound \eqref{eq:bound_integral_points} to \eqref{eq:lattice_points_first_inequality}, we get:
	\begin{equation*} 
			\sum_{v \in \Lambda\setminus \{0\}} r^{- \lVert v \rVert} \leq -\frac{1}{r^\rho} + \frac{\log(r)}{\rho^d} \int_\rho^{+\infty} \frac{(2 t + \rho)^d}{r^t} dt = -\frac{1}{r^\rho} + \frac{2^d r^{\rho/2}}{(\rho \log (r))^d} \int_{3 \rho \log(r)/2}^{+\infty} u^d e^{-u} du,
	\end{equation*}
	where the last equality follows from the change of variables $2u = (2 t + \rho) \log(r)$. Putting $x = 3\rho \log(r)/2$, we recognise the incomplete gamma function \cite[\S~9.2.1]{Erdelyi_1981}:
	\begin{equation*}
	\Gamma(d+1,x) = \int_x^{+\infty} u^d e^{- u} du = d! \, e^{-x} \sum_{k=0}^d \frac{x^k}{k!}.
	\end{equation*}
	The inequality \eqref{eq: bound sum lattice 1} in the Lemma follows. Now under the assumption $x \geq d$, the right hand side of \eqref{eq: bound sum lattice 1} is bounded by
	\begin{equation*}
	\frac{3^d d!}{r^\rho} \sum_{k=0}^d \frac{1}{(d-k)! \, d^k} \leq \frac{3^d d!}{r^\rho} \sum_{k=0}^d \frac{1}{d!} \leq \frac{(d+1) 3^d}{r^\rho}. \qedhere
	\end{equation*}
\end{proof}

To conclude this section, let us see how to deduce the exponential convergence of $m(P_A) \to m(P)$ as $\rho(A) \to +\infty$ from the previous \cref{thm:exponential_convergence_functions}. Note that if a polynomial $P$ does not vanish on the torus $\T^n$, then there is some $\delta >0$ such that it does not vanish on the annulus $\mathcal{C}_\delta$ defined in \eqref{not:C_delta}. 

\begin{corollary} \label{cor:exponential_convergence_polynomials}
	Let $P \in \C[\zn_n^{\pm 1}]$ be a Laurent polynomial that does not vanish on the torus $\T^n$. Then there exist $r>1$ and $C>0$ such that for every matrix $A \in \Z^{m \times n}$ with the property that $\rho(A) \geq 2 d(A)/3\log(r)$, we have the following estimate:
	\begin{equation} \label{eq:exponential_polynomial_bound}
		\lvert m(P_A) - m(P) \rvert \leq \frac{C}{r^{\rho(A)}}.
	\end{equation}
More precisely, one can take $r=e^\delta$, where $\delta > 0$ is any real number such that $P$ does not vanish on the annulus $\mathcal{C}_{\delta}$.
\end{corollary}
\begin{proof}
	Let $P^\ast$ be the conjugate reciprocal of $P$, introduced in \cref{not:polynomials}. Note that for $\zn = (z_1,\dots,z_n) \in \T^n$, we have
	\[
	|P(\zn)|^2 = P(z_1, \ldots, z_n) \overline{P(z_1, \ldots, z_n)} = P(z_1, \ldots, z_n) \overline{P}(\overline{z_1}, \ldots, \overline{z_n}) = P(z_1, \ldots, z_n) \overline{P}\bigl(\frac{1}{z_1}, \ldots, \frac{1}{z_n}\bigr) = PP^\ast(\zn).
	\]
	This shows that $PP^\ast=|P|^2$ on the torus $\T^n$. Fix $\delta>0$ such that $P$ does not vanish on the annulus $\mathcal{C}_{\delta}$. Since the involution $(z_1,\ldots,z_n) \mapsto (\overline{z}_1^{-1},\ldots,\overline{z}_n^{-1})$ preserves $\mathcal{C}_\delta$, the polynomial $P^\ast$ also does not vanish on $\mathcal{C}_\delta$. So, the differential form $\omega := \frac{1}{2} d\log(P P^\ast) = d(P P^\ast)/(2 P P^\ast)$ is holomorphic on some open $U \subseteq \C^n$ containing $\mathcal{C}_\delta$. Moreover, the restriction of $\omega$ to the torus is equal to $d\log |P|$, hence is exact.
	
	Now, for $U \supseteq \mathcal{C}_\delta$ small enough, each loop $\gamma \subseteq U$ is homologous to a loop $\gamma' \subseteq \T^n$. This implies that $\int_\gamma \omega = \int_{\gamma'} d\log\lvert P \rvert =0$. Thus, de Rham's comparison theorem shows that there exists a unique holomorphic function $f \colon U \to \C$ such that $\omega = d f$ on $U$ and $f = \log |P|$ on $\T^n$. 
	Hence $\lvert m(P_A) - m(P) \rvert = \left\lvert \int_{\T^n} f \, d\mu_A - \int_{\T^n} f \, d\mu_n \right\rvert$, and we can apply \cref{thm:exponential_convergence_functions} because $f$ is holomorphic on $U \supseteq \mathcal{C}_\delta$. This yields the bound \eqref{eq:exponential_polynomial_bound}, where we set $C := (n+1) 3^n \cdot \max_{\mathcal{C}_\delta} \lvert f \rvert$.
\end{proof}
\begin{remark}
	For a given Laurent polynomial $P \in \C[\zn_n^{\pm 1}] \setminus \{0\}$ which does not vanish on $\T^n$, one can find an explicit $\delta > 0$, depending on $\min_{\T^n} \lvert P \rvert$, such that $P$ does not vanish on $\mathcal{C}_\delta$. We will carry out this computation for a specific type of polynomial in \cref{pro:holCeps}. 
\end{remark}

\subsection{An explicit error term in the general case}

Let $P  \in \C[\zn_n^{\pm 1}]$ be a Laurent polynomial in $n$ variables, which is not a monomial. Given a matrix $A \in \Z^{m \times n}$, we wish to prove \cref{thm:error-term}, which gives a precise estimate for the error $\lvert m(P_A)-m(P) \rvert$. In order to do so, we approximate the function: 
\[
	\begin{aligned}
	f \colon \T^n &\to \R \cup \{-\infty\} \\
	\zn &\mapsto \log\lvert P(\zn) \rvert
	\end{aligned}
\]
which is singular when $P$ vanishes on $\T^n$, with the smooth functions $f_\varepsilon(\zn) = \frac{1}{2} \log(\lvert P(\zn) \rvert^2 + \varepsilon)$. 
Then, one has that:
\begin{equation} \label{eq:general_bound_decomposition}
	\begin{aligned}
		\lvert m(P_A) - m(P) \rvert \leq &\left| \int_{\T^n} \frac{1}{2} \log(|P|^2+\varepsilon)-\log|P| d\mu_A\right| + \left| \int_{\T^n} \frac{1}{2} \log(|P|^2+\varepsilon)-\log|P| d\mu_n\right| \\ + &\left\lvert \int_{\T^n} f_\varepsilon(\zn) d\mu_A - \int_{\T^n} f_\varepsilon(\zn) d\mu_n(\zn) \right\rvert
	\end{aligned}
\end{equation}
and we proceed by bounding each integral separately. We will show in \cref{pro:holCeps} that $f_\varepsilon$ is holomorphic on a neighbourhood of the torus, so the results of the previous section apply to the last integral. However, two phenomena are competing here. On the one hand, you need to take $\varepsilon$ small enough to make the first two integrals small. On the other hand, the annulus $\mathcal C_\delta$ on which $f_\varepsilon$ is holomorphic becomes smaller and smaller when $\varepsilon \to 0$, which weakens the bound for the third integral given by \cref{thm:exponential_convergence_functions}. 
Thus, the proof of \cref{thm:error-term} will consist in choosing a suitable value of $\varepsilon$, depending on the quantity $\rho(A)$, which balances these two phenomena.

In order to bound the first two integrals, we rely on an explicit estimate for the measure of the set of points $\zn \in \T^n$ where $\lvert P(\zn) \rvert$ is small. 
This estimate is expressed in the following \cref{thm:improved_DH}, which is similar to results by Dobrowolski \cite[Theorem~1.4]{Dobrowolski_2017} and Dimitrov-Habegger \cite[Lemma~A.3]{Dimitrov_Habegger_2019}. 
The main idea of the proof is to compute the aforementioned measure as an integral on $\T^{n-1}$, and then to ``slice'' the torus $\T^{n-1}$ according to the magnitude of a given polynomial $P_\ell(\zn_{n-1})$ related to $P$. This method is due to Habegger, and has already been used in the proof of \cite[Lemma~A.3]{Dimitrov_Habegger_2019}. 
We discuss the relation of our result with other similar results in \cref{rmk:Dimitrov_Habegger}. Finally, we remark that the following \cref{thm:improved_DH} is slightly better than needed for the proof of our main result. As we explain in \cref{rem:sharp}, the order of growth expressed by the bound \eqref{eq:improved_DH} is attained by a specific family of polynomials.

\begin{theorem} \label{thm:improved_DH}
Let $n \geq 1$ be an integer, and let $P \in \C[\zn_n^{\pm 1}]$ be a Laurent polynomial with $k=k(P) \geq 2$ non-zero coefficients. Let moreover $L_\infty(P) := \max_{v \in \Z^n} \lvert c_v(P) \rvert$ and $S(P,r) := \{\zn \in \T^n \colon \lvert P(\zn) \rvert \leq r \}$ for every $r \in \R_{\geq 0}$. Also, set $\kappa = \kappa(P) = \max_{1 \leq i \leq n} (k_i(P))$, where for every $i \in \{1,\dots,n\}$ we let $k_i = k_i(P)$ be the number of non-zero coefficients of $P$ seen as a polynomial in $z_i$. Then for every $\alpha \in ]0,1[$ and every $r>0$, we have the following bound:
\begin{equation} \label{eq:improved_DH}
	\mu_n(S(P,r)) \leq C_1(n,k) \cdot \alpha^{1-n} \Bigl(\frac{r}{L_\infty(P)}\Bigr)^{\frac{1-\alpha}{\kappa-1}},
\end{equation}
where $C_1(n,k) := 6k \cdot (18 n k^2)^{n-1}$.
\end{theorem}
\begin{proof}
	We may assume that $P \in \C[\zn_n]$, since multiplying $P$ by a monomial does not change $S(P,r)$ and $L_\infty(P)$.
	Also, replacing $(P,r)$ by $(P/L_\infty(P),r/L_\infty(P))$, we may assume that $L_\infty(P) = 1$.
	Finally, we may assume $r < 1$, since the right-hand side of \eqref{eq:improved_DH} is $\geq 1$ when $L_\infty(P)=1$ and $r \geq 1$.
	
	We prove the result by induction on $n$. When $n=1$, the result is due to Dobrowolski, as can be seen by applying to \cite[Theorem~1.1]{Dobrowolski_2017} the trivial inequality $(k-1) (\frac{12 \sqrt{2}}{\pi})^{\frac{k-2}{k-1}} \leq C_1(1,k) = 6k$.
	
	If $n \geq 2$, we write $P = \sum_j P_j(\zn_{n-1}) z_n^j$, and we choose $\ell \in \Z$ such that $L_\infty(P_\ell) = L_\infty(P) = 1$. Note that $k(P_\ell) \leq k(P) = k$ and $\kappa(P_\ell) \leq \kappa(P) = \kappa$. Moreover, for every $y \in \T^{n-1}$, we set $Q_y(z) := P(y,z) \in \C[z]$. Similarly, we have $k(Q_y) \leq k_n(P) \leq k$ and $\kappa(Q_y) \leq k_n(P) \leq \kappa$. From now on, we will assume that $P_\ell$ is not a monomial, in other words $\kappa(P_\ell) \geq 2$. We will explain how to cover the case $\kappa(P_\ell)=1$ at the end of the proof.
	
	By Fubini's theorem, we can write:
	\begin{equation} \label{lem max eq 2}
		\mu_n(S(P,r)) = \int_{\T^{n-1}} \mu_1(S(Q_y,r)) d\mu_{n-1}(y) = \int_{|P_\ell(y)| \leq r} \mu_1(S(Q_y,r)) d\mu_{n-1}(y) + \int_{|P_\ell(y)| > r} \mu_1(S(Q_y,r)) d\mu_{n-1}(y) =: I_1 + I_2.
	\end{equation}
	The integral $I_1$ is bounded by $\mu_{n-1}(S(P_\ell,r))$, and by induction we have
	\begin{equation} \label{lem max eq 3}
		I_1 \leq \mu_{n-1}(S(P_\ell,r)) \leq C_1(n-1,k(P_\ell)) \cdot \alpha^{2-n} \Bigl(\frac{r}{L_\infty(P_\ell)}\Bigr)^{\frac{1-\alpha}{\kappa(P_\ell)-1}} \leq C_1(n-1,k) \cdot \alpha^{2-n} r^{\frac{1-\alpha}{\kappa-1}}.
	\end{equation}
	
	We now concentrate on the integral $I_2$. Since $P_\ell(y)$ is a coefficient of $Q_y$, we have $|P_\ell(y)| \leq L_\infty(Q_y)$, so:
	\begin{equation} \label{lem max eq 4}
		\mu_1(S(Q_y,r)) \leq C_1(1,k(Q_y)) \Bigl(\frac{r}{L_\infty(Q_y)}\Bigr)^{\frac{1-\alpha}{\kappa(Q_y)-1}} \leq C_1(1,k) \Bigl(\frac{r}{|P_\ell(y)|} \Bigr)^{\frac{1-\alpha}{\kappa(Q_y)-1}} \leq C_1(1,k) \Bigl(\frac{r}{|P_\ell(y)|} \Bigr)^{\frac{1-\alpha}{\kappa-1}}
	\end{equation}
	for every $y \in \T^{n-1} \setminus S(P_\ell,r)$. Here we assumed $\kappa(Q_y) \geq 2$, but the bound \eqref{lem max eq 4} also holds if $Q_y$ is a monomial, because $L_\infty(Q_y)>r$ implies $S(Q_y,r)=\emptyset$ in this case.
	
	We now slice the domain $\T^{n-1} \setminus S(P_\ell,r)$ according to the modulus of $P_\ell$. Choose an integer $N \geq 1$, and let $r_0,\ldots,r_N$ be the geometric progression starting at $r_0=r$ and ending at $r_N=k(P_\ell)$. Write $\Sigma_m = S(P_\ell,r_m)$ for $m=0,\ldots,N$. Then
	\begin{align*}
		\int_{|P_\ell(y)| > r} \Bigl(\frac{r}{|P_\ell(y)|}\Bigr)^{\frac{1-\alpha}{\kappa-1}} d\mu_{n-1}(y) = \sum_{m=1}^{N} \int_{\Sigma_m \setminus \Sigma_{m-1}} \Bigl(\frac{r}{|P_\ell(y)|}\Bigr)^{\frac{1-\alpha}{\kappa-1}} d\mu_{n-1}(y) & \leq  \sum_{m=1}^{N} \int_{\Sigma_m \setminus \Sigma_{m-1}} \Bigl(\frac{r}{r_{m-1}}\Bigr)^{\frac{1-\alpha}{\kappa-1}} \\
		& \leq \sum_{m=1}^{N} \Bigl(\frac{r}{r_{m-1}}\Bigr)^{\frac{1-\alpha}{\kappa-1}} \cdot (\mu_{n-1}(\Sigma_m) - \mu_{n-1}(\Sigma_{m-1})).
	\end{align*}
	Let $S$ denote the last sum. Using summation by parts, we have
	\begin{equation*}
		S = \Bigl(\frac{r}{r_{N-1}}\Bigr)^{\frac{1-\alpha}{\kappa-1}} \mu_{n-1}(\Sigma_N) - \Bigl(\frac{r}{r_0}\Bigr)^{\frac{1-\alpha}{\kappa-1}} \mu_{n-1}(\Sigma_0) + \sum_{m=1}^{N-1} \mu_{n-1}(\Sigma_m) \cdot \Bigl( \Bigl(\frac{r}{r_{m-1}}\Bigr)^{\frac{1-\alpha}{\kappa-1}} - \Bigl(\frac{r}{r_m}\Bigr)^{\frac{1-\alpha}{\kappa-1}} \Bigr).
	\end{equation*}
	We now bound the volume of $\Sigma_m$ using the induction hypothesis applied to $P_\ell$, with a parameter $\beta \in ]0,\alpha[$ to be specified later:
	\begin{equation*}
		S \leq \Bigl(\frac{r}{r_{N-1}}\Bigr)^{\frac{1-\alpha}{\kappa-1}} + C_1(n-1,k) \beta^{2-n} r^{\frac{1-\alpha}{\kappa-1}} \sum_{m=1}^{N-1} r_m^{\frac{1-\beta}{\kappa(P_\ell)-1}} \bigl( r_{m-1}^{-\frac{1-\alpha}{\kappa-1}} - r_m^{-\frac{1-\alpha}{\kappa-1}} \Bigr).
	\end{equation*}
	Taking the limit when $N$ tends to infinity, the sum over $1 \leq m \leq N-1$ becomes the Riemann--Stieltjes integral of the function $f(x)=x^{(1-\beta)/(\kappa(P_\ell)-1)}$ on the interval $[r,k(P_\ell)]$ with respect to $g(x)=-x^{-(1-\alpha)/(\kappa-1)}$. This integral evaluates to
	\[
	\int_r^{k(P_\ell)} f(x) g'(x) dx = \frac{1-\alpha}{(\kappa-1) \cdot \gamma} \bigl( k(P_\ell)^\gamma - r^\gamma \bigr) \qquad \Bigl(\gamma = \frac{1-\beta}{\kappa(P_\ell)-1} - \frac{1-\alpha}{\kappa-1}\Bigr).
	\]
	Since $\kappa(P_\ell) \leq \kappa$ and $\alpha' < \alpha$, we have $\gamma \geq (\alpha-\beta)/(\kappa-1) > 0$ and also $\gamma < 1$. Using also $1 \leq k(P_\ell) \leq k$, we obtain
	\[
	\int_{|P_\ell(y)| > r} \Bigl(\frac{r}{|P_\ell(y)|}\Bigr)^{\frac{1-\alpha}{\kappa-1}} d\mu_{n-1}(y)
	\leq r^{\frac{1-\alpha}{\kappa-1}} + C_1(n-1,k) \beta^{2-n} \cdot \frac{1-\alpha}{(\kappa-1) \cdot \gamma} k(P_\ell)^\gamma r^{\frac{1-\alpha}{\kappa-1}} \leq r^{\frac{1-\alpha}{\kappa-1}} + k C_1(n-1,k) \frac{\beta^{2-n}}{\alpha-\beta} r^{\frac{1-\alpha}{\kappa-1}}.
	\]
	The function $\beta \in ]0,\alpha[ \mapsto \beta^{2-n} (\alpha-\beta)^{-1}$ has a minimum at $\beta = \frac{n-2}{n-1} \cdot \alpha$ for $n \geq 3$. For this choice of $\beta$, or taking $\beta \to 0$ in the case $n=2$, we get the bound:
	\begin{equation} \label{lem max eq 5}
		\int_{|P_\ell(y)| > r} \Bigl(\frac{r}{|P_\ell(y)|}\Bigr)^{\frac{1-\alpha}{\kappa-1}} d\mu_{n-1}(y)
		\leq \bigl(1+ k C_1(n-1,k) \cdot e \cdot (n-1) \alpha^{1-n}\bigr) r^{\frac{1-\alpha}{\kappa-1}}.
	\end{equation}
	
	Going back to the volume of $S(P,r)$, we obtain using \eqref{lem max eq 2}, \eqref{lem max eq 3}, \eqref{lem max eq 4} and  \eqref{lem max eq 5}:
	\begin{align*}
		\mu_n(S(P,r)) & \leq C_1(n-1,k) \cdot \alpha^{2-n} r^{\frac{1-\alpha}{\kappa-1}} + C_1(1,k) \bigl(1+ k C_1(n-1,k) \cdot e \cdot (n-1) \alpha^{1-n}\bigr) r^{\frac{1-\alpha}{\kappa-1}} \\
		& \leq \bigl( C_1(n-1,k) + C_1(1,k) \cdot e \cdot nk \cdot C_1(n-1,k) \bigr) \alpha^{1-n} r^{\frac{1-\alpha}{\kappa-1}} \\
		& \leq 18 n k^2 \cdot C_1(n-1,k) \alpha^{1-n} r^{\frac{1-\alpha}{\kappa-1}} \\
		& \leq C_1(n,k) \alpha^{1-n} r^{\frac{1-\alpha}{\kappa-1}}.
	\end{align*}
	Finally, let us treat the case $\kappa(P_\ell)=1$. Since $L_\infty(P_\ell)=1$, we have $I_1=0$. Moreover \eqref{lem max eq 4} is still valid, so
	\begin{equation*}
		\mu_n(S(P,r)) = I_2 = \int_{\T^{n-1}} \mu_1(S(Q_y,r)) d \mu_{n-1}(y) \leq C_1(1,k) \cdot r^{\frac{1-\alpha}{\kappa-1}},
	\end{equation*}
	which is less than the right-hand of \eqref{eq:improved_DH}.
\end{proof}
\begin{remark} \label{rmk:Luck}
	We note that Lück \cite[Proposition~2.1]{Luck_2015} provided another estimate for the measure $\mu_n(S(P,t))$. However, this bound depends on the \textit{width} $\mathrm{wd}(P)$, which is a quantity defined in \cite[\S~1.2]{Luck_2015} that turns out to be comparable with $\diam(P)$. In particular $\mathrm{wd}(P_A) \to +\infty$ as $\rho(A) \to +\infty$, which makes Lück's bound not adapted to our purposes.
	More precisely, if we used Lück's bound in the proof of \cref{prop:error_regularization}, we would get an estimate for the first two integrals appearing in \eqref{eq:general_bound_decomposition} which would diverge as $\rho(A) \to +\infty$.
\end{remark}
\begin{remark} \label{rmk:Dimitrov_Habegger}
	We note that Habegger \cite[Lemma~A.4]{Habegger_2018} and Dimitrov and Habegger \cite[Lemma~A.3]{Dimitrov_Habegger_2019} stated another estimate for $\mu_n(S(P,t))$, where the exponent $(1-\alpha)/(\kappa-1)$ appearing in \eqref{eq:improved_DH} is replaced with $1/(2 (k-1))$, but the implicit constant does not depend on $\alpha$.  
	Our method of proof is similar to theirs, and indeed the shape of the bound discussed in \cref{thm:improved_DH} was already suggested after the proof of \cite[Lemma~A.3]{Dimitrov_Habegger_2019}. 
	Finally, Dobrowolski \cite[Theorem~1.4]{Dobrowolski_2017} provided a bound with an exponent of the form $1/\left( \sum_{i = 1}^n (k_i(P) - 1) \right)$, which in general is worse than $1/(\kappa(P) - 1)$.
\end{remark}
\begin{remark}\label{rem:sharp}
	From the bound of \cref{thm:improved_DH}, choosing suitably $\alpha$ we may achieve a bound of the form: 
	\[
	\mu_n(S(P,r)) \ll r^{1/(\kappa-1)} |\log r|^{n-1}
	\] 
	when $r \to 0$. 
	On the other hand, consider the family of polynomials $\{ P_{n,m} := \prod_{i=1}^n (z_i-1)^m \colon n \in \Z_{\geq 1}, \ m\in \Z_{\geq 1} \}$, for which we have that $\kappa(P_{n,m}) = m+1$.
	Then, it can be proved by induction on $n$ and $m$ that the following asymptotic behaviour:
	\[
	\mu_n\left( S\left( P_{n,m},r \right) \right) \asymp  r^{\frac{1}{m}} |\log r|^{n-1}
	\]
	holds true when $r \to 0$. 
	
	Hence, we see that the bound provided by \cref{thm:improved_DH} is actually attained (up to multiplicative constants) by a family of polynomials. However, we expect that for a generic polynomial $P$ it should be possible to get much better upper bounds for the function $\mu_n(S(P,r))$, for instance assuming that the zero locus of $P$ inside $(\C^\times)^n$ intersects $\T^n$ transversally. This will be the subject of future investigations.
\end{remark}

The crucial property of the bound provided by \eqref{eq:improved_DH} is that the constants involved remain bounded if we replace $P$ by $P_A$, for any matrix $A \in \Z^{m \times n}$. 
Under the additional assumptions that $L_\infty(P) = 1$ and $m \leq n$, this suffices to bound the first two integrals appearing in \eqref{eq:general_bound_decomposition}, as we show in the following \cref{prop:error_regularization}. This proposition follows from \cref{thm:improved_DH} by a \textit{Tauberian estimate}, similar in spirit to the ones considered in \cite{Wiener_1932}. 
Note that the aforementioned assumptions are harmless, as we will explain at the beginning of the proof of \cref{thm:error-term}.

\begin{proposition}\label{prop:error_regularization}
	Fix two natural numbers $k \in \Z_{\geq 2}$ and $n \in \Z_{\geq 1}$, and a real number $0<\alpha<1$. Let $P \in \C[\zn_n^{\pm 1}]$ be a Laurent polynomial such that $k(P) = k$ and $L_\infty(P) = 1$. 
	Let $A$ be a matrix in $\Z^{m \times n}$ such that $\rho(A) > \diam(P)$ and $m \leq n$. Then, for every $\varepsilon>0$ the following inequalities hold:
	\[ 
	0 \leq \int_{\T^n} \frac{1}{2} \log(|P|^2+\varepsilon)-\log|P| d\mu_A \leq  C_2(n,k,\alpha) \cdot \varepsilon^{\frac{1-\alpha}{2(k-1)}}
	\]
	where $C_2(n,k,\alpha) := 12k^2(18 n k^2)^{n-1} \cdot \frac{\alpha^{1-n}}{1-\alpha}$.
\end{proposition}

\begin{proof}
	Let $\nu$ be the measure on $\R_{\geq 0}$ defined as the push-forward of $\mu_A$ along the (measurable) function $\lvert P \rvert \colon \T^n \to \R_{\geq 0}$.
	Moreover, let us define the functions $\phi(t) := \frac{1}{2} \log\left( 1 + \frac{\varepsilon}{t^2} \right)$ and $\psi(s) := \sqrt{\frac{\varepsilon}{e^{2s} - 1}}$, so that $\psi(\phi(x)) = \phi(\psi(x)) = x$ for every $x > 0$. Writing $\chi_S$ for the characteristic function of a subset $S \subseteq \R$, we have the following identities:
	\begin{equation} \label{eq:Tauberian_equalities}
		\begin{aligned}
			\int_{\T^n} \frac{1}{2}\log(|P|^2+\varepsilon)-\log|P| d\mu_A &= \int_0^{+\infty} \phi(t) d \nu(t) = \int_0^{+\infty} \left( \int_0^{+\infty} \chi_{[0,\phi(t)]}(s) \, ds \right) d\nu(t) \\
			&= \int_0^{+\infty} \left( \int_0^{+\infty} \chi_{[0,\psi(s)]}(t) \, d\nu(t) \right) ds = \int_0^{+\infty} \nu([0,\psi(s)]) ds.
		\end{aligned}
	\end{equation}
	Now, we can bound $\nu([0,t])$ using \cref{thm:improved_DH}. Indeed, the assumption $\rho(A) > \diam(P)$ implies that $k(P_A) = k(P) = k$ and $L_\infty(P_A) = L_\infty(P) = 1$, as we explained at the beginning of the proof of \cref{cor:uniform L2}.
	Hence, we see from \cref{thm:improved_DH} that:
	\[
	\nu([0,t]) := \mu_A(\{ \zn \in \T^n \colon \lvert P(\zn) \rvert \leq t \}) = \mu_m(\{ \zn \in \T^m \colon \lvert P_A(\zn) \rvert \leq t \}) \leq C_1\left(m,k(P_A)\right) \cdot \alpha^{1-m} \cdot t^{\frac{1-\alpha}{\kappa(P_A)-1}} \leq C_1(n,k) \cdot \alpha^{1-n} \cdot t^{\frac{1-\alpha}{k-1}}
	\]
	because $\kappa(P_A) \leq k(P_A) = k$ and $C_1$ is an increasing function in both variables. Combining this bound with the identities provided by \eqref{eq:Tauberian_equalities}, we get:
	\[
		\int_{\T^n} \frac{1}{2}\log(|P|^2+\varepsilon)-\log|P| d\mu_A = \int_0^{+\infty} \nu([0,\psi(s)]) ds \leq \left( C_1(n,k) \cdot \alpha^{1-n} \cdot \int_0^{+\infty} \left(e^{2 s} - 1\right)^{-\frac{1-\alpha}{2  (k-1)}} d s \right) \cdot \varepsilon^{\frac{1-\alpha}{2 (k-1)}}.
	\]
	The last integral can be computed by substituting $u = e^{2 s} - 1$ and using Cauchy's residue theorem (see \cite[p.~107]{Cartan_1961}):
	\[
		\int_0^{+\infty} \left(e^{2 s} - 1\right)^{-\frac{1-\alpha}{2 (k-1)}} d s = \frac{1}{2} \int_0^{+\infty} \frac{u^{-\frac{1-\alpha}{2  (k-1)}}}{u+1} du = \frac{\pi}{2 \cdot \sin\left( \frac{\pi(1-\alpha)}{2  (k-1)} \right)}.
	\]
To conclude, we observe that $2 \sin(x) \geq x$ when $0 \leq x \leq \pi/2$ and that $2 \frac{k-1}{1-\alpha}\cdot C_1(n,k) \cdot \alpha^{1-n} \leq C_2(n,k,\alpha)$.
\end{proof}

We now tackle the last term in \eqref{eq:general_bound_decomposition}. We wish to bound it using \cref{thm:exponential_convergence_functions}. To this end, we must show that the function $\frac{1}{2} \log(\lvert P \rvert^2 + \varepsilon)$, defined on the torus $\T^n$, can be extended to a holomorphic function on a neighbourhood of $\T^n$. Following the idea used to prove \cref{cor:exponential_convergence_polynomials}, we use the Laurent polynomial $PP^\ast$, which is equal to $|P|^2$ on $\T^n$, and we consider the function $\frac12 \log(PP^\ast+\varepsilon)$. The following proposition shows that this function is well-defined and holomorphic on an explicit annulus containing $\T^n$.

\begin{proposition} \label{pro:holCeps}
	Fix $n \in \N$, a Laurent polynomial $P \in \C[\zn_n^{\pm 1}] \setminus \{0\}$, and $\varepsilon>0$. Let $P^\ast$ be the conjugate reciprocal of $P$, introduced in \cref{not:polynomials}. Let $\delta := \delta_\varepsilon(P)$ be the constant defined in \eqref{not:a_eps}.
	Then the function $f_\varepsilon := \frac{1}{2} \log(P P^\ast + \varepsilon)$, defined using the principal branch of the logarithm, is holomorphic on an open neighbourhood of the annulus $\mathcal{C}_\delta$ defined in \eqref{not:C_delta}.
\end{proposition}
\begin{proof}
	Denote by $Q$ the Laurent polynomial $PP^\ast$. Since $Q=|P|^2$ on the torus, the image $Q(\T^n)$ is a segment contained in $\R_{\geq 0}$. We will show that $Q(\mathcal{C}_\delta)$ is contained in the half-plane $\{w \in \C : \mathrm{Re}(w)>-\varepsilon\}$, so that we may consider the principal branch of the logarithm of $Q+\varepsilon$ on $\mathcal{C}_\delta$. Our strategy is to fix a point $\underline{u} = (u_1,\ldots,u_n)$ on $\T^n$, and to bound from below the real part of $Q$ near $\underline{u}$ in radial directions. We therefore write $\zn = (e^{h_1} u_1, \ldots, e^{h_n} u_n)$, and consider the function $g_{\underline{u}} : \R^n \to \C$ defined by
	\begin{equation*}
	g_{\underline{u}}(\underline{h}) = Q(e^{h_1} u_1, \ldots, e^{h_n} u_n).
	\end{equation*}
	Note in particular that $g_{\underline{u}}(\mathbf{0}) = Q(\underline{u}) \in \R_{\geq 0}$.
	Let us apply the multivariable Taylor theorem to $g_{\underline{u}}$ at the origin:
	\begin{equation*}
		g_{\underline{u}}(\underline{h}) = g_{\underline{u}}(\mathbf{0}) + \sum_{i=1}^n \frac{\partial g_{\underline{u}}}{\partial h_i}(\mathbf{0}) \cdot h_i + R_{\underline{u}}(\underline{h})
	\end{equation*}
where the remainder is given in Lagrange’s form by
\begin{equation*}
R_{\underline{u}}(\underline{h}) = \frac12 \sum_{i,j=1}^n \frac{\partial^2 g_{\underline{u}}}{\partial h_i \partial h_j}(\alpha_{\underline{h}} \, \underline{h}) \cdot h_i h_j \qquad (0 < \alpha_{\underline{h}} <1).
\end{equation*}

\begin{lemma} \label{lem:Dg_imag}
The partial derivatives $\partial g_{\underline{u}}/\partial h_i(\mathbf{0})$, $1 \leq i \leq n$, are purely imaginary.
\end{lemma}

\begin{proof}
From the definition of $Q$, we have $\overline{Q(z_1,\ldots,z_n)} = Q(1/\overline{z_1}, \ldots,1/\overline{z_n})$ for every $\zn \in (\C^\times)^n$. Substituting $\zn = (e^{h_1} u_1,\dots,e^{h_n} u_n)$ gives $\overline{g_{\underline{u}}(\underline{h})} = g_{\underline{u}}(-\underline{h})$. Differentiating with respect to $h_i$ at $\mathbf{0}$, we get the result.
\end{proof}

\cref{lem:Dg_imag} ensures that the real part of $Q$ behaves \emph{quadratically} in $\underline{h}$ near the torus. More precisely, we have
	\begin{equation} \label{eq:holomorphicity_bound_1}
		\mathrm{Re}(Q(\zn)) = \mathrm{Re}(g_{\underline{u}}(\underline{h})) = g_{\underline{u}}(\mathbf{0}) + \mathrm{Re}(R_{\underline{u}}(\underline{h})) \geq \mathrm{Re}(R_{\underline{u}}(\underline{h})) \geq - \lvert R_{\underline{u}}(\underline{h}) \rvert.
	\end{equation}
	Now, let us introduce the differential operators $D_k = z_k (\partial/\partial z_k)$ for every $k \in \{1,\dots,n\}$. Notice that $\partial g_{\underline{u}}/\partial h_k (\underline{h}) = (D_k Q)(\zn)$, and similarly for the higher order derivatives. So, in order to give an upper bound for $\lvert R_{\underline{u}}(\underline{h}) \rvert$, it suffices to bound $\lvert D_i D_j Q(\underline{z}) \rvert$. Expanding this polynomial, we get
	\begin{equation*}
	D_i D_j Q(\zn) = D_i D_j \biggl(\sum_{v, {w} \in \Z^n} c_{v}(P) \overline{c_{w}(P)} \zn^{v - {w}}\biggr) = \sum_{v, {w} \in \Z^n} (v_i-w_i) (v_j-w_j) c_v(P) \overline{c_{w}(P)} \zn^{v - {w}}
	\end{equation*}
	which gives the bound
	\begin{equation} \label{eq:holomorphicity_bound_2}
	\lvert D_i D_j Q(\zn) \rvert \leq \sum_{v, {w} \in \Z^n} |v_i-w_i| |v_j-w_j| |c_v(P)| |c_w(P)| \prod_{k=1}^n e^{|h_k| |v_k-w_k|} \leq \diam(P)^2 L_1(P)^2 e^{\diam(P) \lVert \underline{h} \rVert_1}.
	\end{equation}
	Combining \eqref{eq:holomorphicity_bound_1} and \eqref{eq:holomorphicity_bound_2}, we obtain
	\begin{equation} \label{eq:holomorphicity_bound_3}
		\mathrm{Re}(Q(\zn)) \geq - \lvert R_{\underline{u}}(\underline{h}) \rvert \geq - \frac{1}{2} \sum_{i,j=1}^n \left\lvert \frac{\partial^2 g_{\underline{u}}}{\partial h_i \partial h_j}(\alpha_{\underline{h}} \, \underline{h}) \right\rvert \cdot \lvert h_i \rvert \lvert h_j \rvert \geq - \frac12 \diam(P)^2 L_1(P)^2 e^{\diam(P) \lVert \underline{h} \rVert_1} \cdot \lVert \underline{h} \rVert_1^2.
	\end{equation}
	Finally, if $\zn \in \mathcal{C}_{\delta}$ then $\lVert \underline{h} \rVert_1 = \sum_{i = 1}^n \lvert \log\lvert z_i \rvert \rvert \leq \delta$, and the definition \eqref{not:a_eps} of $\delta$ implies that the right-hand side of \eqref{eq:holomorphicity_bound_3} is $\geq -\frac{2}{3}\varepsilon$. We thus have $\mathrm{Re}(Q+\varepsilon)>0$ on an open neighbourhood of $\mathcal{C}_{\delta}$, as we wanted to show.
\end{proof}

Thanks to \cref{pro:holCeps}, we may apply \cref{thm:exponential_convergence_functions} to the functions $f_\varepsilon$. However, we also need to bound $f_\varepsilon$ on the domain $\mathcal{C}_{\delta_\varepsilon(P)}$. This is the content of the following lemma.

\begin{lemma} \label{lem:estimate_value_feps}
	Let $P \in \C[\zn_n^{\pm 1}]$ be a non-zero Laurent polynomial, and fix $\varepsilon > 0$. Let $\delta := \delta_\varepsilon(P) \in \R_{> 0}$ be defined as in \eqref{not:a_eps}, and let $f_\varepsilon$ be the function defined in \cref{pro:holCeps}. Then, we have that: 
	\[
		\lvert f_\varepsilon(\zn) \rvert \leq \lvert \log \varepsilon \rvert + 2 \lvert \log L_1(P) \rvert + 3
	\] 
	for every $\zn \in \mathcal{C}_\delta$.
\end{lemma}
\begin{proof}
	Let $Q = P P^\ast$. We have $f_\varepsilon = \log(Q+\varepsilon)=\log |Q+\varepsilon| + i \, \mathrm{arg}(Q+\varepsilon)$, and by the proof of \cref{pro:holCeps}, the argument of $Q+\varepsilon$ stays in $]-\pi/2,\pi/2[$ on the domain $\mathcal{C}_\delta$. It remains to bound from below and from above the modulus of $Q+\varepsilon$.
	
	The lower bound follows from \eqref{eq:holomorphicity_bound_3}, since $|Q(\zn)+\varepsilon| \geq \mathrm{Re}(Q(\zn))+\varepsilon \geq \varepsilon/3$ for $\zn \in \mathcal C_\delta$ as seen at the end of the previous proof. For the upper bound, let us write $\zn = (e^{h_1} u_1,\dots,e^{h_n} u_n)$, where $\underline{u} = (u_1,\dots,u_n) \in \T^n$ and $\underline{h} = (h_1,\dots,h_n) \in \R^n$. Then, a simple application of the triangle inequality yields:
	\begin{equation*}
	\lvert Q(\zn) \rvert \leq \sum_{v, w \in \Z^n} \lvert c_v(P) \rvert \lvert c_w(P) \rvert \prod_{j=1}^n e^{\lvert h_j \rvert \lvert v_j - w_j \rvert} \leq L_1(P)^2 e^{\delta_\varepsilon(P) \diam(P)} \leq \frac{4}{3} L_1(P)^2.
	\end{equation*}
	We deduce that:
	\begin{equation*}
	\lvert \log(Q(\zn)+\varepsilon) \rvert \leq \lvert \log \lvert Q(\zn)+\varepsilon\rvert \rvert + \frac{\pi}{2} \leq \max\left(\, \left\lvert \log \frac{\varepsilon}{3} \right\rvert, \left\lvert \log\left(\frac{4}{3} L_1(P)^2+\varepsilon\right) \right\rvert \, \right) + \frac{\pi}{2}.
	\end{equation*}
	We conclude using the inequality $|\log(x+y)| \leq |\log(x)| + |\log(y)| + \log 2$, valid for any $x,y>0$, which is easily proved by distinguishing the cases $x \leq y$ and $x \geq y$.
\end{proof}

Using \cref{thm:exponential_convergence_functions}, together with \cref{pro:holCeps} and the bound of \cref{lem:estimate_value_feps}, we get:
\begin{corollary}\label{coro:bound2}
Let $m,n \geq 1$ be integers, and $P(\zn_n) \in \C[\zn_n^{\pm 1}]$ be a non-zero Laurent polynomial. For every $\varepsilon > 0$, consider the function $f_\varepsilon = \frac 12\log(|P|^2+\varepsilon)$, and let $\delta_\varepsilon(P) \in \R_{> 0}$ be defined as in \eqref{not:a_eps}. Then, for every matrix $A \in \Z^{m \times n}$ such that the inequality $\rho(A) \delta_\varepsilon(P) \geq 2 d(A)/3$ holds true, we have:
\[
\left|\int_{\T^n} f_\varepsilon d\mu_A - \int_{\T^n} f_\varepsilon d\mu_n\right| \leq (d(A)+1) 3^{d(A)} \cdot \frac{|\log \varepsilon| + 2 |\log L_1(P)| + 3}{\exp(\delta_\varepsilon(P) \cdot \rho(A))}. 
\]
\end{corollary}

We are finally ready to prove \cref{thm:error-term}, by choosing a suitable value of $\epsilon$.
\begin{proof}[Proof of \cref{thm:error-term}]
We fix the Laurent polynomial $P$ and a matrix $A$ in $\Z^{m\times n}$ verifying the assumption of \cref{thm:error-term}, that is $\rho(A)\geq \rho_0(P)$. For this proof, we will simplify a bit the notations, by denoting $k:=k(P)$, $\rho_0:=\rho_0(P)$, $\delta_\epsilon:=\delta_\epsilon(P)$, $\rho := \rho(A)$ and $d = d(A)$. Note that $d\leq n$ by definition.

Since the quantity $\lvert m(P)-m(P_A) \rvert$, which we want to bound, does not change when multiplying $P$ by a non-zero constant, we will assume without loss of generality that $L_\infty(P) = 1$.
Finally, we may assume without loss of generality that $m \leq n$ (compare with \cite[Theorem~4]{Smyth_2018}). 
Indeed, let $U \in \Z^{m \times m}$ be a matrix with non-zero determinant such that $U \cdot A$ is in row-echelon form, and write $B \in \Z^{m' \times n}$ for the matrix obtained from $U \cdot A$ by deleting all the zero rows. 
Then, $m' = \mathrm{rk}(B) \leq n$ by construction, and one has that $\Lambda_A = \Lambda_B$, which implies that $\rho(A) = \rho(B)$. Moreover, one sees that $m(P_A) = m((P_{A})_U) = m(P_{U \cdot A}) = m(P_B)$ by combining \cite[Lemma~7]{Smyth_2002} and \cite[Lemma~6]{Smyth_2018}.
Hence, upon replacing $A$ with $B$, we can and will assume without loss of generality that $m \leq n$, as we claimed.

Now, let us come back to the bound provided by \eqref{eq:general_bound_decomposition}:
\begin{equation*}
	\lvert m(P_A) - m(P) \rvert \leq \left| \int_{\T^n} \frac{1}{2} \log(|P|^2+\varepsilon)-\log|P| d\mu_A\right| + \left| \int_{\T^n} \frac{1}{2} \log(|P|^2+\varepsilon)-\log|P| d\mu_n\right| + \left\lvert \int_{\T^n} f_\varepsilon(\zn) d\mu_A - \int_{\T^n} f_\varepsilon(\zn) d\mu_n(\zn) \right\rvert
\end{equation*}
which holds for any $\epsilon>0$. 
The first two terms are bounded by \cref{prop:error_regularization} under some conditions, whereas the third term is bounded by \cref{coro:bound2}, under other assumptions. The strategy is to choose the value of $\epsilon$ such that we can indeed apply these two results and, moreover, that the two upper bounds become comparable.
To that end, let us fix the quantity: 
\begin{equation} \label{eq:our_epsilon}
	\varepsilon := \left( \frac{(1-\alpha)\diam(P) \cdot L_1(P)}{k-1} \cdot \frac{\log(\rho)}{\rho} \right)^2
\end{equation}
where $0 \leq \alpha \leq \frac{1}{2}$ is a real number to be fixed later.

To begin with, we may apply \cref{prop:error_regularization}, since $\rho\geq \rho_0 > \diam(P)$ by assumption (see \eqref{not:rho_0} for the definition of $\rho_0$), as well as $L_\infty(P)=1$ and $m \leq n$.

Next, to verify the assumptions of \cref{coro:bound2}, we have to check that $\delta_\varepsilon \rho \geq 2 d / 3$. 
To see this, note first of all that $\rho\geq \rho_0\geq 7\mathrm{diam}(P)^2$.
Plugging this bound in \eqref{eq:our_epsilon}, and using the elementary inequalities $\frac{\log(\rho)}{\rho}\leq \frac{3}{4\sqrt{\rho}}$ and $k-1 >1-\alpha$, we see that $\sqrt{\epsilon} \leq \log(4/3) L_1(P)$. Combining this with the definition of $\delta_\epsilon$, given in \eqref{not:a_eps}, we get the equality:
\begin{equation} \label{eq:delta epsilon rho}
\delta_\epsilon = \frac{\sqrt{\varepsilon}}{\mathrm{diam}(P) \cdot L_1(P)} = \frac{1-\alpha}{k-1} \cdot \frac{\log(\rho)}{\rho}. 
\end{equation}
The desired condition $\delta_\epsilon \rho \geq 2d/3$ then follows from the lower bound $\rho \geq \rho_0 \geq \exp(2 (k-1) \max(n,5)) \geq \exp\left(\frac{k-1}{1-\alpha}d\right)$, as we have both $1-\alpha \geq \frac 12$ and $d \leq n$.

Eventually, to use efficiently \cref{coro:bound2}, we need an upper bound on the quantity $\lvert \log(\epsilon) \rvert +2 \log(L_1(P)) +3$. We begin by noting that our assumption that $\rho\geq \rho_0 \geq 7 \diam(P)^2$, together with the remark $1 = L_\infty(P) \leq L_1(P) \leq k L_\infty(P) = k$ and the elementary inequality $\frac{\log(\rho)}{\rho} \leq \frac{3}{4 \sqrt{\rho}}$, implies that:
\begin{equation*}
\epsilon \leq  \left( \frac{(1-\alpha)\diam(P) \cdot L_1(P)}{k-1}\right)^2  \cdot \frac{9}{16\rho} \leq \left( \frac{k}{k-1} \right)^2 \cdot \frac{9}{112} \leq 1.
\end{equation*}
This allows us to bound $|\log \varepsilon|$. Using the expression \eqref{eq:delta epsilon rho} for $\varepsilon$ and our assumption $\log(\rho) \geq 2 (k-1) \max(n,5) \geq \frac{ e^{3/2} \cdot (k-1)}{(1-\alpha)\diam(P)}$, we get:
\[
	|\log(\varepsilon) \rvert + 2 \log(L_1(P)) + 3 = 2 \log\biggl( \frac{e^{3/2} (k-1)}{(1-\alpha)\diam(P)} \cdot \frac{\rho}{\log(\rho)} \biggr) \leq 2 \log(\rho).
\] 

Applying \cref{prop:error_regularization} and \cref{coro:bound2} together with the last estimate, we obtain the following bound:
\begin{equation*} 
	\lvert m(P_A) - m(P) \rvert \leq C_2(n,k,\alpha) \cdot \varepsilon^{\frac{1-\alpha}{2(k-1)}} + (d + 1) 3^{d} \cdot \frac{2 \log(\rho)}{\exp(\delta_\varepsilon \rho)}.
\end{equation*}
We can now use our definition \eqref{eq:our_epsilon} of $\epsilon$, together with the bounds $L_1(P) \leq k$ and $d\leq n$, to obtain:
\begin{equation*}
	\begin{aligned}
	\lvert m(P_A) - m(P) \rvert &\leq C_2(n,k,\alpha) \cdot \left(\frac{(1-\alpha)\diam(P) \cdot L_1(P)}{k-1} \cdot \frac{\log(\rho)}{\rho}\right)^{\frac{1-\alpha}{k-1}} + (d + 1) 3^d \cdot \frac{2\log(\rho)}{\rho^{\frac{1-\alpha}{k-1}}} \\
	&\leq \left(C_2(n,k,\alpha) \cdot \left(\frac{\diam(P) \cdot L_1(P)}{k-1}\right)^{\frac{1-\alpha}{k-1}}+ 2(n + 1) 3^n \right)\frac{\log(\rho)}{\rho^\frac{1-\alpha}{k-1}} \\ 
	&\leq \left(C_2(n,k,\alpha) \cdot \left(\frac{k}{k-1}\right)^{\frac{1}{k-1}}+ 2(n + 1) 3^n \right) \cdot \log(\rho) \cdot \left( \frac{\diam(P)}{\rho} \right)^\frac{1-\alpha}{k-1}.
	\end{aligned}
\end{equation*}
The constant of the last inequality can be bounded using elementary inequalities, giving
\begin{equation} \label{eq:final_bound}
	\lvert m(P_A) - m(P) \rvert \leq 25 k^2 (18nk^2)^{n-1} \cdot \frac{\alpha^{1-n}}{1-\alpha} \cdot \log(\rho) \cdot \left( \frac{\diam(P)}{\rho} \right)^\frac{1-\alpha}{k-1}.
\end{equation}
Finally, we convert \eqref{eq:final_bound} into a bound with a logarithmic factor $\log(\rho)^n$. To this end, we observe that the function: 
\[
\alpha \mapsto \alpha^{1-n} \cdot \rho^{\frac{\alpha}{k-1}}
\] 
achieves its minimum at the value $\alpha = \frac{n (k-1)}{\log(\rho)}$. Note that $0<\alpha \leq \frac 12$, because $\rho\geq \rho_0\geq e^{2n(k-1)}$, and $\rho^{\frac{\alpha}{k-1}} = e^n$. Thus, we obtain:
\[
	\begin{aligned}
	\lvert m(P_A) - m(P) \rvert
	& \leq 25 k^2 (18nk^2)^{n-1} \cdot \frac{ \log(\rho)^{n-1}}{(1-\alpha) (n(k-1))^{n-1}} \cdot \log(\rho) \cdot e^n \cdot \left( \frac{\diam(P)}{\rho} \right)^\frac{1}{k-1} \\
	& \leq \frac{50 k^2 (18nk^2)^{n-1} e^n}{(n(k-1))^{n-1}} \cdot \log(\rho)^n \cdot \left( \frac{\diam(P)}{\rho} \right)^\frac{1}{k-1} \\
	& \leq 8 \cdot (36ek)^{n-1} \cdot \log(\rho)^n \cdot \left( \frac{\diam(P)}{\rho} \right)^\frac{1}{k-1}.
	\end{aligned}
\]
This gives the theorem.
\end{proof}

\begin{remark}
As is clear from the proof of \cref{thm:error-term}, the quality of the error term depends essentially only on the quality of the bound \eqref{eq:improved_DH} provided by \cref{thm:improved_DH}. To see this, fix a Laurent polynomial $P \in \C[\zn_n^{\pm 1}]$ with $k(P) \geq 2$, and a matrix $A \in \Z^{m \times n}$ such that $m \leq n$ and $\rho(A) > \diam(P)+1$, so that $k(P) = k(P_A)$ and $L_\infty(P) = L_\infty(P_A)$. 
Suppose moreover that there exist $a, b, C \in \R_{> 0}$ such that the bounds: 
\begin{equation*} 
	\mu_n(S(P,t)) \leq C \cdot \alpha^{-b} \cdot (t/L_\infty(P))^{(1-\alpha) \cdot a} \quad \text{and} \quad \mu_m(S(P_A,t)) \leq C \cdot \alpha^{-b} \cdot (t/L_\infty(P_A))^{(1-\alpha) \cdot a}
\end{equation*}
hold true for every $t \in \R_{> 0}$ and $0 < \alpha < 1$.
Then, going through the proof of \cref{thm:error-term} one sees that there exist two constants $C',\rho_0 \geq 0$ such that the following bound holds: 
\begin{equation*} 
	\lvert m(P_A) - m(P) \rvert \leq C' \cdot \log(\rho(A))^{b+1} \cdot \left(\frac{\diam(P)}{\rho(A)}\right)^a
\end{equation*}
for every matrix $A$ such that $\rho(A) \geq \rho_0$.
\end{remark}

\section{Discussion of the speed of convergence}
\label{sec:speed_of_convergence}
	
In this section, we study several situations where more can be said about the error term $m(P_A)-m(P)$, compared to the bounds given in \cref{cor:exponential_convergence_polynomials} and \cref{thm:error-term}. 
In particular, we devote \cref{sec:Condon} to an experimental study of the differences $m(P(z_1,z_1^d)) - m(P)$ for some two-variable polynomials $P \in \Z[\zn_2]$. A full asymptotic expansion for these sequences, under a technical assumption on $P$, was provided by Condon \cite{Condon_2012}, and our experiments are compatible with this result. Finally, we devote \cref{sec:Pd} to the study  of a multivariate example where not only an equivalent of the error term can be obtained, but also a full asymptotic expansion.
This example goes beyond Condon's framework, both with respect to the number of variables involved (as we study a family of $2$-variable polynomials whose Mahler measures converge to the Mahler measure of a $4$-variable one) and the type of expansion that we get, where a logarithmic term appears.	
	
\subsection{Asymptotic expansions in the presence of toric points}
\label{sec:Condon}

\begin{figure}[t]
	\centering
	\begin{subfigure}[b]{0.32\linewidth}
		\centering
		\includegraphics[width=\linewidth]{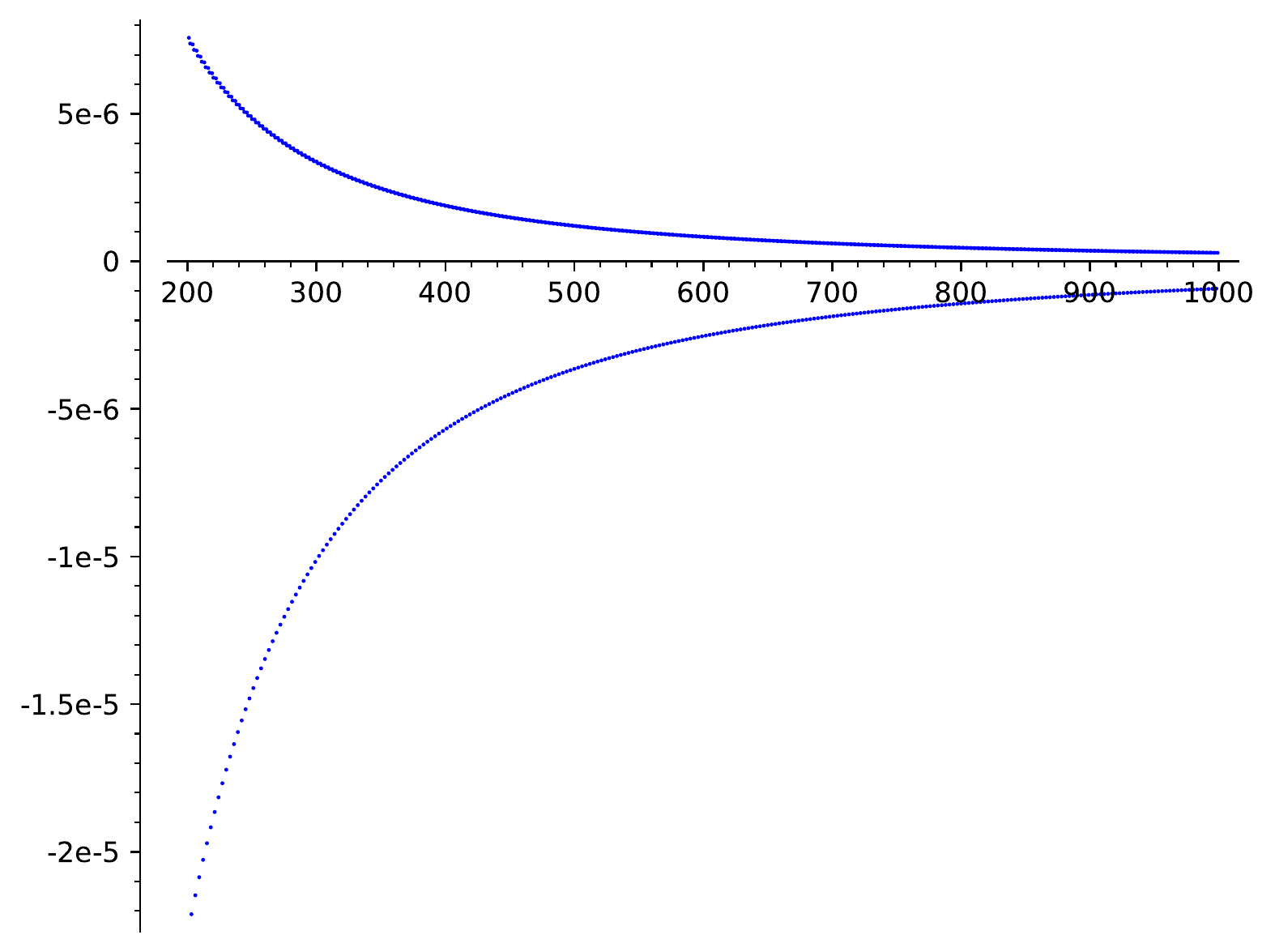}
		\caption{{\small $ P(z_1,z_2) = z_2 + (z_1 + 1) $}}
		\label{fig:Dirichlet_A}
	\end{subfigure}
	\begin{subfigure}[b]{0.32\linewidth}
		\centering
		\includegraphics[width=\linewidth]{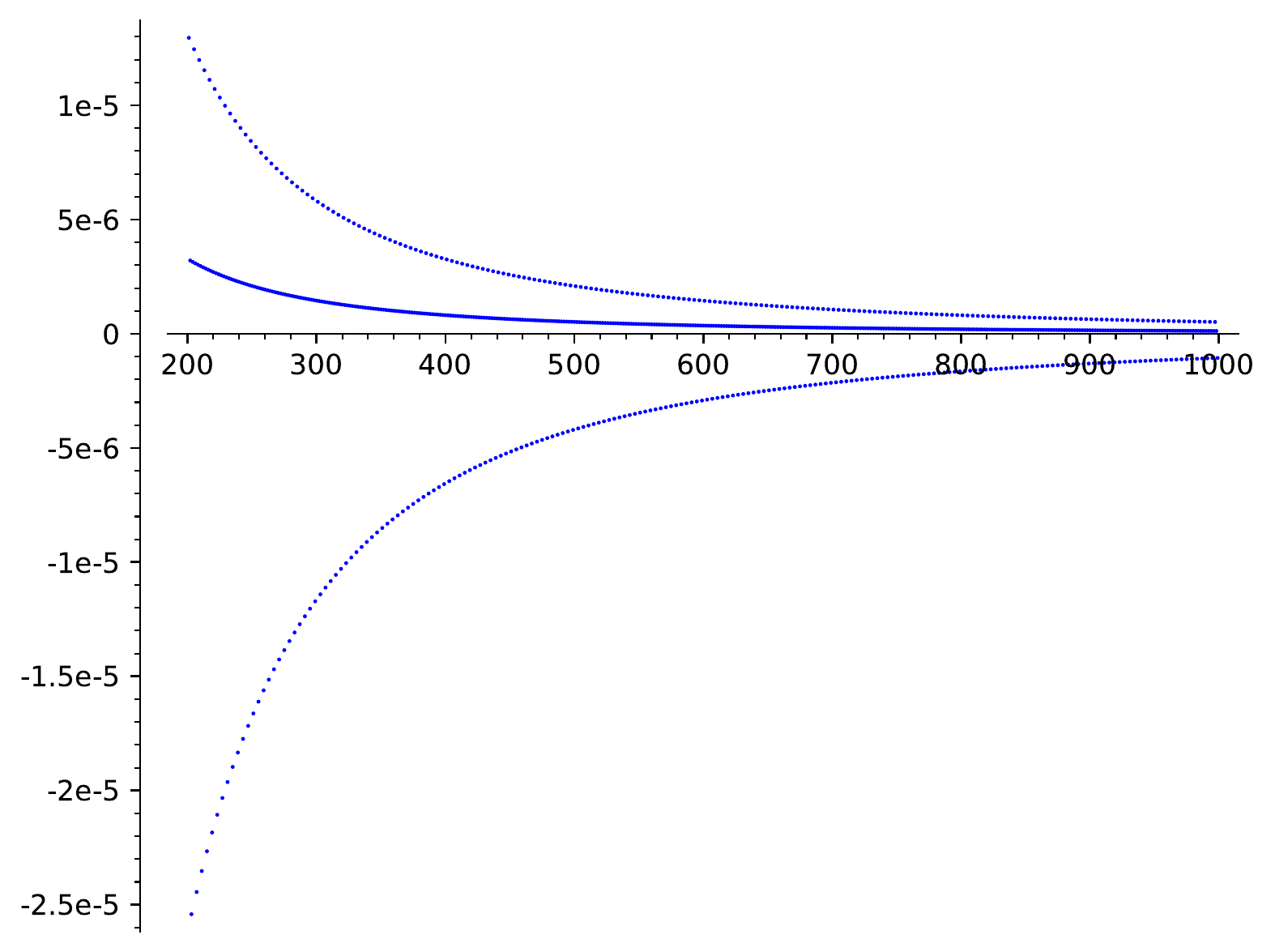}
		\caption{{\small $ P(z_1,z_2) =  (z_1 + 1) z_2 + (z_1 - 1) $}}
		\label{fig:Dirichlet_B}
	\end{subfigure}
	\begin{subfigure}[b]{0.32\linewidth}
		\centering
		\includegraphics[width=\linewidth]{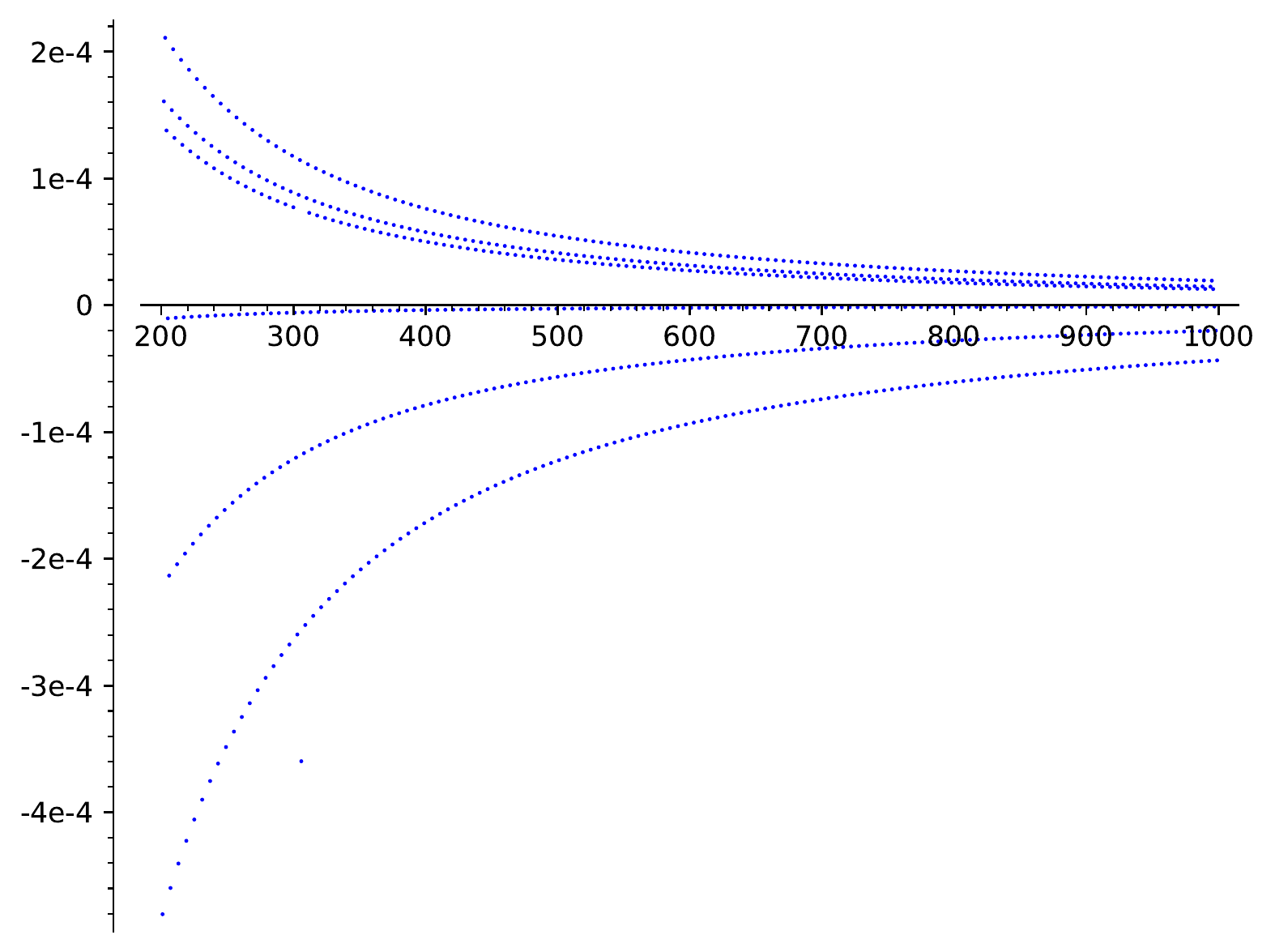}
		\caption{{\small $P(z_1,z_2) = z_1 z_2^2 + (z_1^2 + z_1 + 1) z_2 + z_1$}}
		\label{fig:Deninger}
	\end{subfigure}
	
	\captionsetup{justification=centering}
	\caption{Plots of $m(P_{A_d}) - m(P)$, for $A_d = (1,d)$ and $200 \leq d \leq 1000$, which seem to lie on finitely many smooth branches. \newline 
	The \textsc{SageMath} code we used to produce these plots is available online \cite{Guilloux_2021}.}
	\label{fig:periodic_examples}
\end{figure}

When a polynomial $P(\zn_n) \in \C[\zn_n^{\pm 1}]$ vanishes on the torus $\T^n$, one cannot hope that the error term $\lvert m(P_A) - m(P) \rvert$ decays exponentially fast as $\rho(A) \to +\infty$.
This is already evident when $n = 2$, and we take the sequence of matrices $A_d := (1,d) \in \Z^{1 \times 2}$, which results in the sequence of polynomials $P_{A_d}(z_1) := P(z_1,z_1^d)$. For $P(z_1,z_2) = z_1 + z_2 + 1$, the resulting sequence of polynomials $\{ P_{A_d}(z_1) = z_1 + z_1^d + 1 \}_{d = 1}^{+\infty}$ was already studied by Boyd \cite[Appendix~2]{Boyd_1981}, who proved that:
\[
	m(z_1 + z_1^d + 1) - m(z_1 + z_2 + 1) = \frac{c(d)}{d^2} + O\left( \frac{1}{d^3} \right)
\]
where $c \colon \Z \to \R$ is a $3$-periodic function.
More precisely, $c(d) := - \sqrt{3} \pi/6$ if $d \equiv 2 (3)$, and $c(d) := \sqrt{3} \pi/18$ otherwise. This is reflected by the fact that the plot of $m(P_{A_d}) - m(P)$, depicted in \cref{fig:Dirichlet_A}, consists of two branches.
This is by no means an isolated phenomenon: we include in \cref{fig:periodic_examples} two other examples, taken from \cite[Equation~(1-7)]{Boyd_1998} and \cite[Table~1]{Boyd_1998} respectively, of polynomials $P$ for which the error term $m(P_{A_d}) - m(P)$ appears to be divided into a finite number of smooth branches.

\begin{figure}[b]
	\centering
	\begin{subfigure}[t]{0.32\linewidth}
		\centering
		\includegraphics[width=\linewidth]{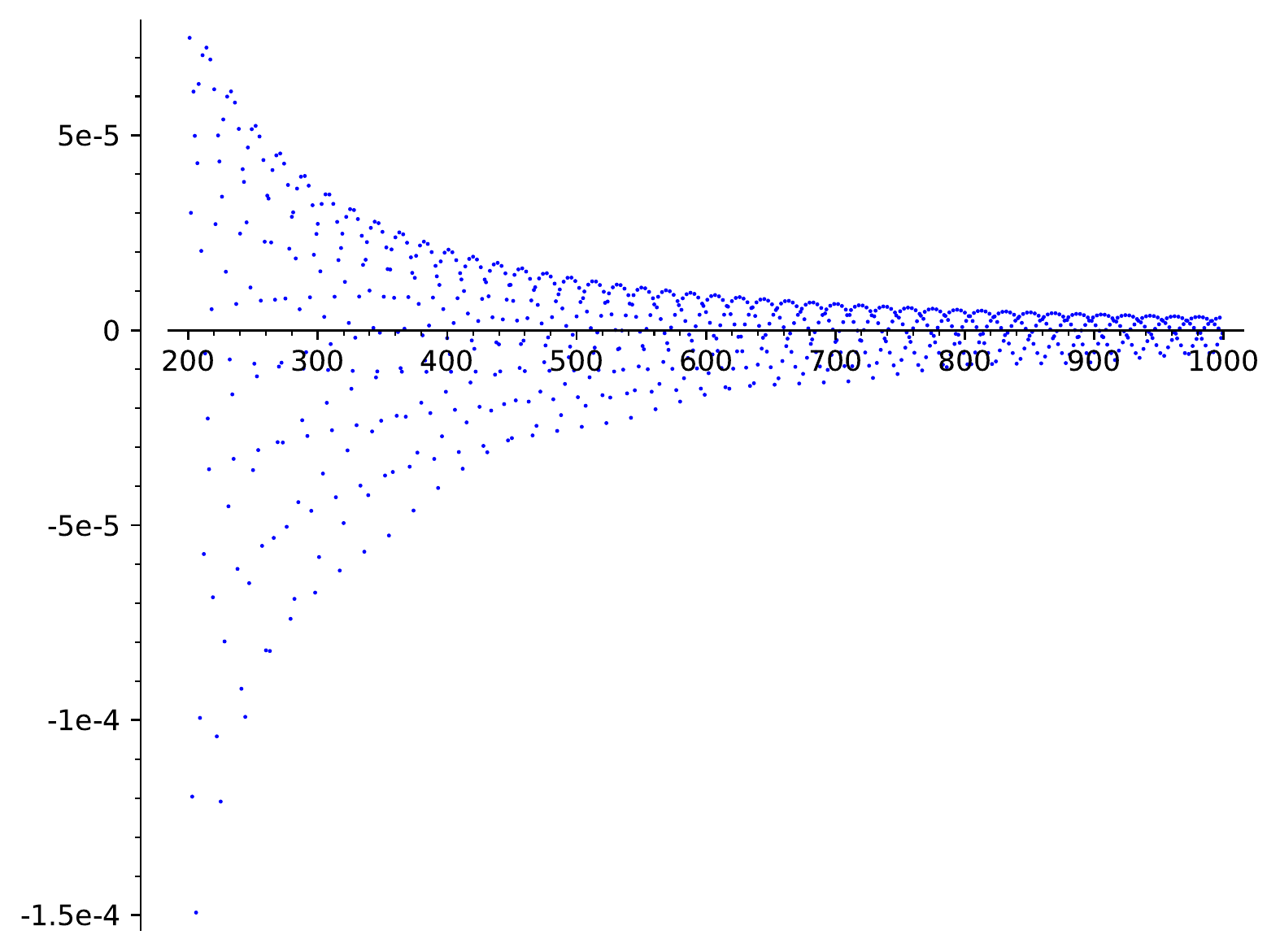}
		\caption{{\small $P = (z_1+1)^4 z_2 - (z_1^2+1) (z_1^2-z_1+1) $}}
		\label{fig:Dirichlet_C}
	\end{subfigure} 
	\begin{subfigure}[t]{0.32\linewidth}
		\centering
		\includegraphics[width=\linewidth]{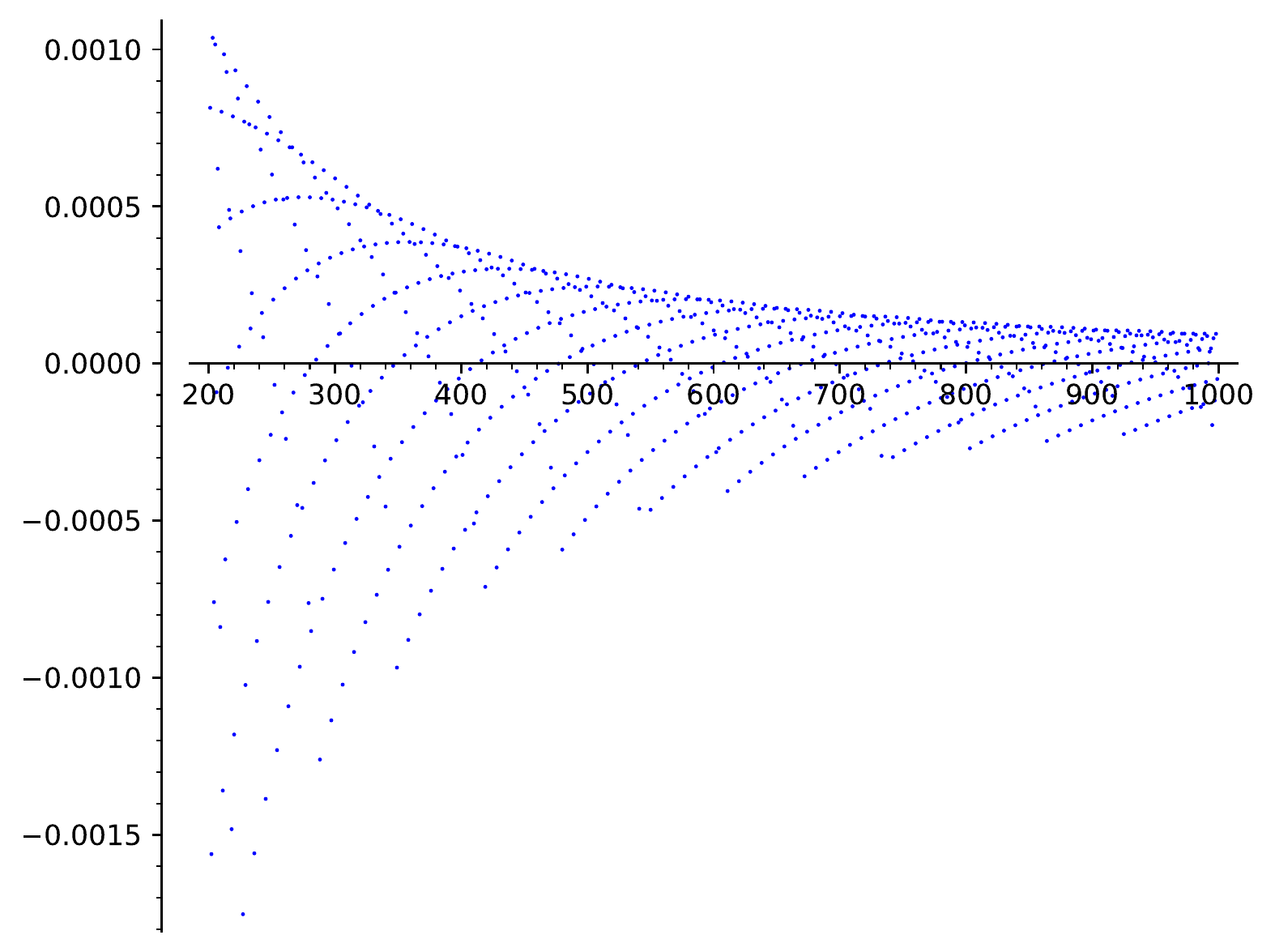}
		\caption{{\small $ P = (z_1^8 + z_1^6 + z_1^4 + z_1^2 + 1) (z_2^2 + 1) + \newline (2 z_1^8 - 37 z_1^6 + 5 z_1^5 + 70 z_1^4 + 5 z_1^3 - 37 z_1^2 + 2) z_2$}}
		\label{fig:Dirichlet_D}
	\end{subfigure}
	\begin{subfigure}[t]{0.32\linewidth}
		\centering
		\includegraphics[width=\linewidth]{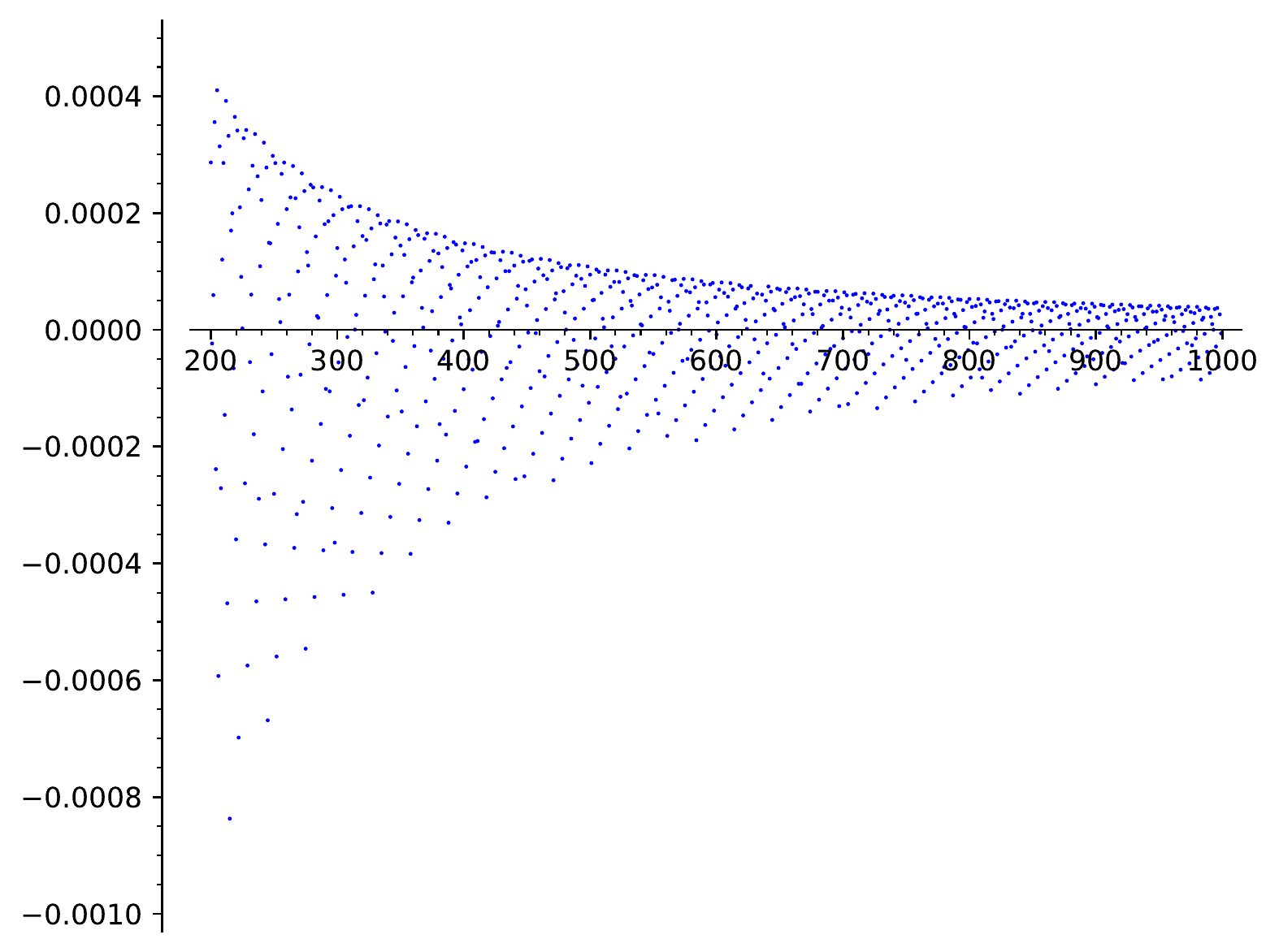}
		\caption{{\small $P = (z_1 + 1) z_2^2 + (z_1^2 + z_1 + 1) z_2 + z_1^2 + z_1$}}
		\label{fig:Mellit}
	\end{subfigure}
	
	\captionsetup{justification=centering}
	\caption{Plots of $m(P_{A_d}) - m(P)$, for $A_d = (1,d)$ and $200 \leq d \leq 1000$, which seem to lie on infinitely many smooth branches. \newline 
		The \textsc{SageMath} code we used to produce these plots is available online \cite{Guilloux_2021}.}
	\label{fig:quasi-periodic_examples}
\end{figure}

Note however that not all polynomials $P(z_1,z_2)$ give rise to an error term $m(P_{A_d}) - m(P)$ with this kind of behaviour. This is depicted in \cref{fig:Dirichlet_C,fig:Dirichlet_D,fig:Mellit}, which display polynomials taken from \cite[Table~1]{Boyd_Rodriguez-Villegas_2001}, \cite[Example~4.8]{Liu_Qin_2021} and \cite[Equation~1]{Mellit_2019} respectively.

These different types of phenomena have been partially explained by Condon's work \cite{Condon_2012}, which provides an asymptotic expansion for the error term $m(P_{A_d}) - m(P)$ of an irreducible polynomial $P \in \C[\zn_2^{\pm 1}]$ such that $P$ and $\partial P/\partial z_2$ do not have a common root on $\T^2$.
To be more precise, we need to recall some terminology introduced by \cite{Condon_2012}. First of all, a function $c \colon \R \to \R$ is said to be \textit{quasi-periodic} if it is the sum of finitely many continuous periodic functions. Then, for every collection of quasi-periodic functions $\{ c_j \colon \R \to \R \}_{j \in \N}$ and every function $f \colon \N \to \R$ we will use the notation $f(d) \approx \sum_{j \in \N} \frac{c_j(d)}{d^j}$ if, for every $J \in \N$, there exists a constant $C_{f,J} > 0$ such that the following bound:
\[
	\left\lvert f(d) - \sum_{j = 0}^{J} \frac{c_j(d)}{d^j} \right\rvert \leq \frac{C_{f,J}}{d^{J+1}}
\]
holds for every $d \geq 1$. This notation generalises the usual notion of asymptotic series (see for example \cite[Definition~1.3.1]{Bleistein_Handelsman_1986}), where the coefficients $c_j$ are assumed to be constant.
In particular, \cite[Proposition~1]{Condon_2012} shows that, as in the classical case, any given function $f \colon \N \to \R$ has at most one asymptotic expansion of this kind.
Then, Condon proves in \cite[Theorem~1]{Condon_2012} that for any Laurent polynomial $P \in \C[\zn_2^{\pm 1}]$ such that $P$ and $\partial P/\partial z_2$ do not have common zeros on $\T^2$, one has an asymptotic expansion:
\begin{equation} \label{eq:Condon}
	m(P(z_1,z_1^d)) - m(P) \approx \sum_{j=2}^{+\infty} \frac{c_j(d)}{d^j}
\end{equation}
where each $c_j \colon \R \to \R$ is an explicit quasi-periodic function, given by a linear combination of the periodic functions: 
\[
\{ t \mapsto \mathfrak{B}_k( \langle \theta - t \varphi \rangle ) \colon k \in \{2,\dots,j\}, \ (e^{2 \pi i \theta},e^{2 \pi i \varphi}) \in V_P(\C) \cap \T^2 \}
\] 
where $\mathfrak{B}_k(x)$ denotes the $k$-th Bernoulli polynomial, and $\langle x \rangle := x - \lfloor x \rfloor$ denotes the fractional part of a real number $x \in \R$.
In particular, if $V_P(\C) \cap \T^2 \subseteq \mu_N \times \mu_N$, where $\mu_N \subseteq \T^1$ denotes the set of $N$-th roots of unity, each function $c_j$ is $N$-periodic. This is precisely what happens for the two polynomials displayed in \cref{fig:Dirichlet_A,fig:Dirichlet_B}, and hence this periodicity of the coefficients $c_j$ explains why each point $(d,m(P_{A_d}) - m(P)) \in \R^2$ seems to lie on a finite union of graphs of smooth functions.
On the other hand, one can show that each of the polynomials displayed in \cref{fig:quasi-periodic_examples} has toric points whose coordinates are not roots of unity, and this gives rise to the depicted behaviour, where the points $(d,m(P_{A_d}) - m(P)) \in \R^2$ seem to lie on an infinite union of graphs of smooth functions.
Note finally that \cref{fig:Deninger,fig:Dirichlet_D,fig:Mellit} do not fall strictly within the framework of \cite[Theorem~1]{Condon_2012}, because in these cases $P$ and $\partial P/\partial z_2$ have common roots on $\T^2$.
This is related to the fact that $m(P_{A_d}) - m(P)$ seems to decay more slowly than $1/d^2$ in some cases. For instance, extensive computational evidence (already mentioned in \cite[\S~8.3]{Condon_2012}), shows that for the polynomial $P$ appearing in \cref{fig:Deninger} one might expect that $m(P_{A_d}) - m(P) \sim c(d)/d^{3/2}$, where $c \colon \Z \to \R$ is $6$-periodic.	
	
\subsection{An asymptotic expansion with a logarithmic term}
\label{sec:Pd}
	
\begin{figure}[b]
\centering
\includegraphics[width=.4\linewidth]{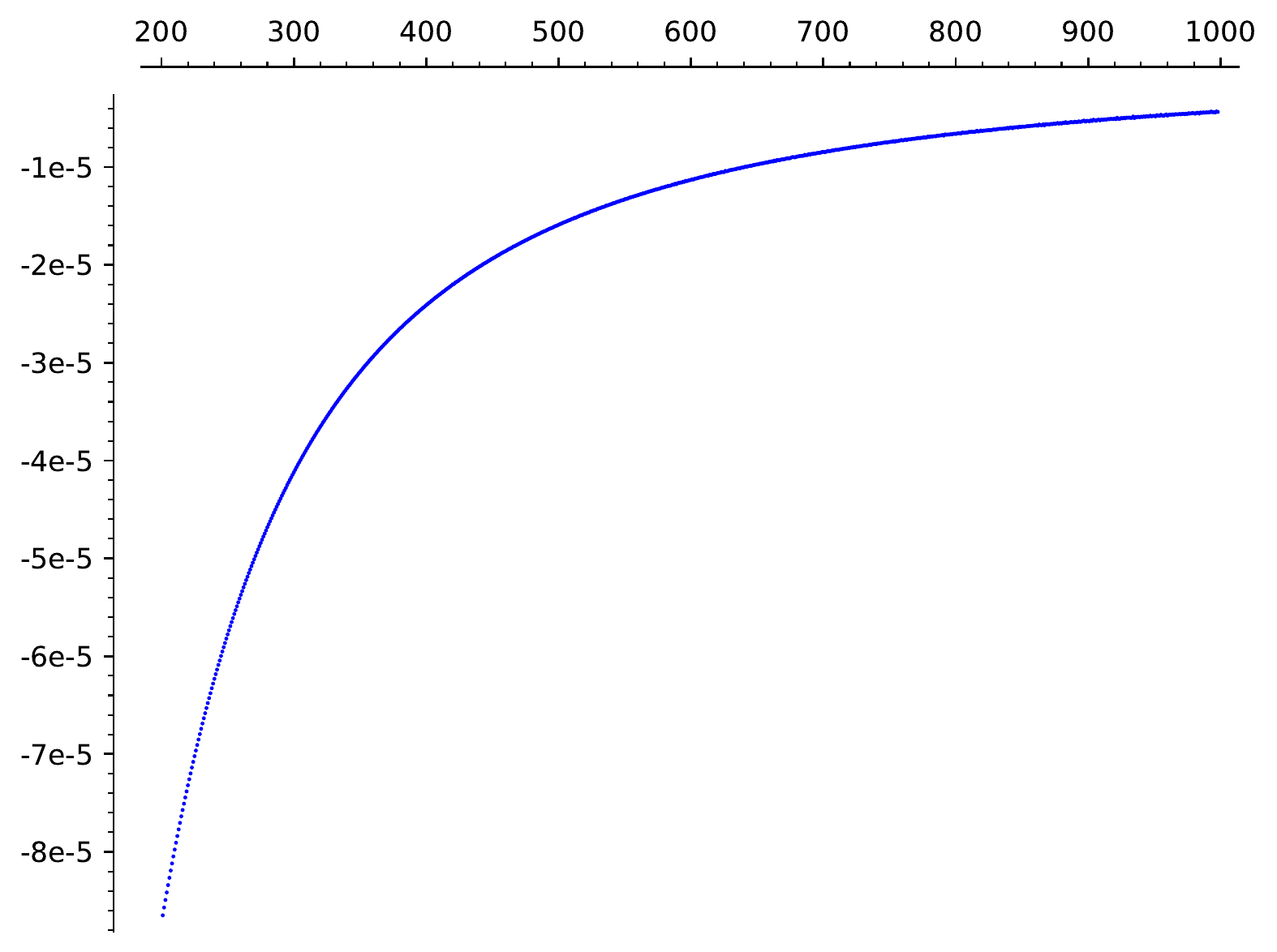}
\captionsetup{justification=centering,margin=5cm}
\caption{Plot of $m(P_d) - m(P_\infty)$, for $200 \leq d \leq 1000$. \newline
	The \textsc{SageMath} code we used to produce this plot is available online \cite{Guilloux_2021}.}\label{fig:mPd}
\end{figure}

This section is dedicated to the sequence of polynomials $ P_d(z_1,z_2) := \sum_{0 \leq i + j \leq d} z_1^i z_2^j \in \C[z_1,z_2] $, whose Mahler measure was widely studied in \cite{Mehrabdollahei_2021} by the third author of this paper. 
In particular, she proved that
\begin{equation} \label{eq:limit_P_d}
		\lim_{d \to +\infty} m(P_d) = \frac{9}{2 \pi^2} \zeta(3) = -18 \cdot \zeta'(-2)
\end{equation}
	where $ \zeta(s) $ denotes Riemann's zeta function. This convergence is illustrated in \cref{fig:mPd} and exhibits a much simpler behaviour than the examples discussed before.
	
	We can give a new proof of \eqref{eq:limit_P_d} using \cref{thm:Lawton_generalisation}.
	More precisely, we can write:
	\[
		P_d(\zn_2) = \frac{1}{(1 - z_1) (1 - z_2)} - \left( \frac{z_1}{(1 - z_1) (z_1 - z_2)} \right) z_1^{d+1} - \left( \frac{z_2}{(1 - z_2) (z_2 - z_1)} \right) z_2^{d+1}
	\]
	using the geometric series. Thus, we see that:
	\[
		P_d(\zn_2) (1 - z_1) (1 - z_2) (z_1 - z_2) = z_1^{d + 2} (z_2 - 1) + z_2^{d + 2} (1 - z_1) + (z_1 - z_2) = (1 - z_1^{d+2}) (1  - z_2) - (1 - z_1) (1 - z_2^{d+2})
	\]
	which implies that $m(P_d) = m(P_\infty(\zn_2^{M_d}))$, where $P_\infty(\zn_4) := (1-z_1)(1-z_2) - (1-z_3)(1-z_4)$ and $M_d := \begin{psmallmatrix}
	d + 2 & 0 & 1 & 0 \\ 0 & 1 & 0 & d + 2
	\end{psmallmatrix} \in \Z^{2 \times 4}$.
	To apply \cref{thm:Lawton_generalisation}, we need to compute $\rho(M_d)$. This is elementary:
	\[
		\{ v \in \Z^4 \mid M_d \cdot v = \mathbf{0} \} = \left\langle \begin{psmallmatrix}
		-1 \\ 0 \\ d + 2 \\ 0
		\end{psmallmatrix}, \begin{psmallmatrix}
		0 \\ d + 2 \\ 0 \\ -1
		\end{psmallmatrix} \right\rangle_\Z
	\]
	so that $ \rho(M_d) = d + 2 $ for every $ d \in \N $.
	Thus, \cref{thm:Lawton_generalisation} shows that: 
	\[
		\lim_{d \to +\infty} m(P_d) = \lim_{d \to +\infty} m(P_\infty(\zn_2^{M_d})) = m(P_\infty).
	\]
	Finally, D'Andrea and Lalín \cite[Theorem~7]{D'Andrea_Lalin_2007} have proved that $m(P_\infty) = -18 \cdot \zeta'(-2)$, which yields back the convergence \eqref{eq:limit_P_d} proved by Mehrabdollahei in \cite{Mehrabdollahei_2021}.
	
	The proof of \eqref{eq:limit_P_d} provided in \cite{Mehrabdollahei_2021} proceeds along very different lines. 
	More precisely, Mehrabdollahei uses crucially the fact that $P_d$ is always an exact polynomial (see \cite[Definition~2.2]{Mehrabdollahei_2021}), which allows her to write:
	\begin{equation} \label{eq:mPd_formula}
		m(P_d) = \frac{3}{d+1} \sum_{1 \leq k \leq d+1} \frac{(d+2-2k)}{2\pi} \cdot D\left(e^{\frac{2\pi ik}{d+2}}\right) - \frac{3}{d+2} \sum_{1 \leq k \leq d} \frac{(d+1-2k)}{2\pi} \cdot D\left(e^{\frac{2\pi ik}{d+1}}\right)
	\end{equation}
	where $D(z) := \arg(1-z) \log\lvert z \rvert - \mathrm{Im}\left( \int_0^z \log(1-t) \frac{dt}{t} \right)$ denotes the Bloch-Wigner dilogarithm (see \cite[Theorem p.~2]{Mehrabdollahei_2021}).
	Even if our proof of \eqref{eq:limit_P_d} does not use \eqref{eq:mPd_formula}, the latter identity allows us to obtain a full asymptotic expansion for the error term $m(P_d) - m(P_\infty)$, which is the content of the following theorem.
	
\begin{theorem} \label{thm:Pd_asymptotics}
	We have the following asymptotic expansion of $m(P_d) - m(P_\infty)$ as $d \to +\infty$:
	\begin{equation} \label{eq:mPd_asymptotics}
		 m(P_d) - m(P_\infty) \approx \frac{1}{(d+1) (d+2)} \left[ -\frac{\log(d)}{2} + \sum_{k=0}^{+\infty} \frac{\alpha_k}{d^k} \right]
	\end{equation}
	where the coefficients $\alpha_k \in \R$ are defined as:
	\[
		\begin{aligned}
			\alpha_0 &:= 6 \bigl(\zeta'(-1) - \zeta'(-2)\bigr) + \frac{\log(2 \pi)}{2} - 1 \\
			\alpha_k &:= \frac{12 \cdot (-1)^k}{k (k+1)} \sum_{j = 0}^{\lfloor k/2 \rfloor} \binom{k+1}{2 j} \cdot \frac{(2^{k+1- 2j} - 1) (2 j-1)}{(2j+1)(2j+2)} \cdot B_{2j+2} \cdot \zeta(2 j) \qquad (k \geq 1)
		\end{aligned}
	\]
	where $B_n$ denotes the $n$-th Bernoulli number. 
	
	In particular, we have that $\displaystyle m(P_d)-m(P_\infty) \sim -\frac{\log(\rho(M_d))}{2\rho(M_d)^2}$ as $d \to +\infty$. 
\end{theorem}
    \begin{proof}
    	First of all, we observe that \eqref{eq:mPd_formula} can be rewritten as:
    	\begin{equation} \label{eq:mPd_error}
    		m(P_d) - m(P_\infty) = \frac{3}{(d+1) (d+2)} \left[ -2 \zeta'(-2) + (d+1)^3 E_{d+1}(f) - (d+2)^3 E_{d+2}(f) \right]
    	\end{equation}
    	where $f(x) := \frac{1 - 2 x}{2 \pi} D(e^{2 \pi i x})$ and $E_d(F) := \int_0^1 F(x) d x - \frac{1}{d} \left( \sum_{j = 1}^{d} F\left(\frac{j}{d}\right) \right)$ for each integrable function $F \colon [0,1] \to \R$.
	    	
    	Now, write $h(x) := (2 x - 1) x \log(x)$ and observe that the function 
    	\[
    		g(x) := f(x) - \log(2 \pi) (1 - 2 x)^2  - h(x) - h(1-x)
    	\] 
    	is smooth on the closed interval $[0,1]$. 
    	Indeed, $f(x)$ and $g(x)$ are smooth on the open interval $(0,1)$, and $g(x) = g(1-x)$. Thus, to see that $g(x)$ is smooth on $[0,1]$ it is sufficient to compute the Maclaurin series:
    	\begin{equation} \label{eq:Maclaurin_g}
    		g(x) = (1 - 2 x) \left[ -\log(2 \pi) + (\log(2 \pi) + 2) x + \sum_{k = 2}^{+\infty} \frac{1}{(1 - k) k} \cdot x^k + \sum_{m = 1}^{+\infty} \frac{\zeta(2 m)}{m (2 m + 1)} \cdot x^{2 m + 1} \right]
    	\end{equation}
    	which follows from the identity $\frac{d^2}{d \theta^2} D(e^{2 i \theta}) = -2 \cot(\theta)$. Moreover, \eqref{eq:Maclaurin_g} allows one to write the asymptotic expansion:
    	\[
    		E_d(g) \approx \frac{3 \log(2 \pi) + 2}{6} \cdot \frac{1}{d^2} + \sum_{m = 2}^{+\infty} \frac{B_{2 m} [2 m + 1 + 2 (2 m - 3) \zeta(2 m - 2)]}{m (2 m - 1) (2 m - 2) (2 m - 3)} \cdot \frac{1}{d^{2 m}}
    	\]
    	using the classical Euler-Maclaurin summation formula \cite[Equation~1]{Navot_1962}.
	    	
    	This formula was extended by Navot to functions with a logarithmic singularity at one endpoint of the integration interval \cite[Equation~7]{Navot_1962}.
    	Applying this generalisation to $h(x)$ we see that:
    	\[
    		E_d(h) \approx \frac{1}{12} \cdot \frac{\log(d)}{d^2} - \left( \zeta'(-1) + \frac{1}{12} \right) \cdot \frac{1}{d^2} + 2 \zeta'(-2) \cdot \frac{1}{d^3} - \sum_{m = 2}^{+\infty} \frac{B_{2 m} (2 m + 1)}{(2 m) (2 m - 1) (2 m - 2) (2 m - 3)} \cdot \frac{1}{d^{2 m}}
    	\]
    	as follows from the Taylor expansion $h(x) = (x-1) + \frac{5}{2} (x - 1)^2 + \sum_{k = 3}^{+\infty} \frac{(-1)^{k+1} (k+2)}{k (k-1) (k-2)} \cdot (x-1)^k$.
    	
    	Hence, observing that $E_d((1 - 2 x)^2) = -\frac{2}{3} \cdot \frac{1}{d^2}$ and $E_d(h(x)) = E_d(h(1-x))$, we get:
    	\begin{equation} \label{eq:error_f}
    		E_d(f) \approx \frac{1}{6} \cdot \frac{\log(d)}{d^2} + \left( \frac{1 - \log(2 \pi)}{6} - 2 \zeta'(-1) \right) \cdot \frac{1}{d^2} + 4 \zeta'(-2) \cdot \frac{1}{d^3} + 4 \sum_{m = 2}^{+\infty} a_{2 m} \cdot \frac{1}{d^{2 m}}
    	\end{equation}
    	where we set $a_k := \frac{B_k \zeta(k - 2)}{k (k - 1) (k - 2)} \in \Q \cdot \pi^{k - 2}$ for every integer $k \geq 4$.
		
	Now, combining the identities:
    	\begin{align}
\nonumber	    		(d+1) \log(d+1) - (d+2) \log(d+2) &\approx -\log(d) -1 + \sum_{k=1}^{+\infty} \frac{(-1)^k (2^{k+1} - 1)}{k (k+1)} \cdot \frac{1}{d^k} \\
\nonumber   		\frac{1}{(d+1)^{2 m - 3}} - \frac{1}{(d+2)^{2 m - 3}} &\approx (2 m - 3) \sum_{j = 1}^{+\infty} \binom{j + 2m - 4}{2m - 3} \frac{(-1)^{j + 1} (2^j - 1)}{j} \cdot \frac{1}{d^{j + 2 m - 3}} 
    	\end{align}
    	with \eqref{eq:mPd_error} and \eqref{eq:error_f} we get:
    	\[
\begin{aligned}
    		(d+1) (d+2) (m(P_d) - m(P_\infty)) &\approx -\frac{\log(d)}{2} + \left( 6 (\zeta'(-1) - \zeta'(-2)) + \frac{\log(2 \pi)}{2} - 1 \right) + \\ 
			&+ \sum_{k=1}^{+\infty} \frac{(-1)^k (2^{k+1} - 1)}{2 k (k+1)} \cdot \frac{1}{d^k} \\&+ 12 \sum_{k=2}^{+\infty} (-1)^k \cdot \left( \sum_{j = 1}^{k - 1} \binom{k-1}{k-j} \frac{(2^j - 1)(k-j)}{j} a_{k-j+3} \right) \frac{1}{d^k}
\end{aligned}
    	\]
which after some rearrangement, gives us \eqref{eq:mPd_asymptotics}.
\end{proof}
\begin{remark}
	The asymptotic expansion \eqref{eq:mPd_asymptotics} has been checked numerically using the PARI/GP program \texttt{Asympraw} available at \cite{Belabas_Cohen_2021}.
\end{remark}
\begin{remark}
	We note that in the asymptotic expansion \eqref{eq:mPd_asymptotics}, the coefficients $\alpha_k$ do not depend on $d$, which is in contrast to what happened for the examples described in \cref{sec:Condon}. 
	Moreover, if $k \geq 1$ we see that $\alpha_k$ is a $\Q$-linear combination of $1, \pi^2, \pi^4, \ldots, \pi^{2 \lfloor k/2 \rfloor}$, whereas $\alpha_0$ and $\pi$ are most likely algebraically independent.
\end{remark}

\begin{remark} \label{rmk:Gu_Lalin}
	Note that \cite[Proposition~8]{Gu_Lalin_2021} provides another family of polynomials, in three variables, whose Mahler measures converge to $m(P_\infty)$. They correspond to the monomial substitutions provided by the matrices:
	\[
		A_{a,b} := \begin{pmatrix}
		 b & 0 & 0 & a \\
		 0 & 1 & 0 & 0 \\
		 0 & 0 & 1 & 0
		\end{pmatrix}
		\in \Z^{3 \times 4}
	\]
	taken as either $a \to +\infty$ or $b \to +\infty$, where $a,b \in \N$ are coprime.
	Since $\ker(A_{a,b}) \cap \Z^4 = \Z \cdot (-a,0,0,b)^t$, we see that $\rho(A_{a,b}) = \max(a,b)$, and so \cite[Proposition~8]{Gu_Lalin_2021} can be seen as a special case of \cref{thm:Lawton_generalisation}. On the other hand, the proof provided by Gu and Lalín uses an explicit formula (see \cite[Theorem~1]{Gu_Lalin_2021}) for the Mahler measures of the three-variable polynomials $(P_\infty)_{A_{a,b}}$, which is similar to the formula \eqref{eq:mPd_formula} proved in \cite{Mehrabdollahei_2021} by the third named author of this paper.
\end{remark}

\subsection{Perspectives}

We hope that the previous \cref{sec:Condon,sec:Pd} managed to convey to the reader our impression that understanding the rate of convergence (and, even more, the asymptotic expansions) of the difference $m(P_A) - m(P)$, remains a difficult and interesting challenge. 
In particular, the bound provided by \cref{thm:error-term} seems far from optimal, even for a general polynomial.
Moreover, the actual rate of convergence, for a fixed polynomial $P$, seems to depend on the geometry of the real algebraic set $V_P(\C) \cap \T^n$, which can be quite complicated on its own (see \cite[Example~5.2.5]{Guilloux_Marche_2021}). Furthermore, one should study as well the geometries of the intersections of this real algebraic set with the sub-tori cut out by the matrices $A$.
For example, it would be interesting to explain the exponent $3/2$ observed in the asymptotics for $P$ given in \cref{fig:Deninger}. We also lack a rationale explaining the logarithmic term appearing in the asymptotic expansion provided in \cref{thm:Pd_asymptotics}.

\begin{figure}[b]
	\centering
	\begin{subfigure}[t]{0.45\linewidth}
		\centering
		\includegraphics[width=\linewidth]{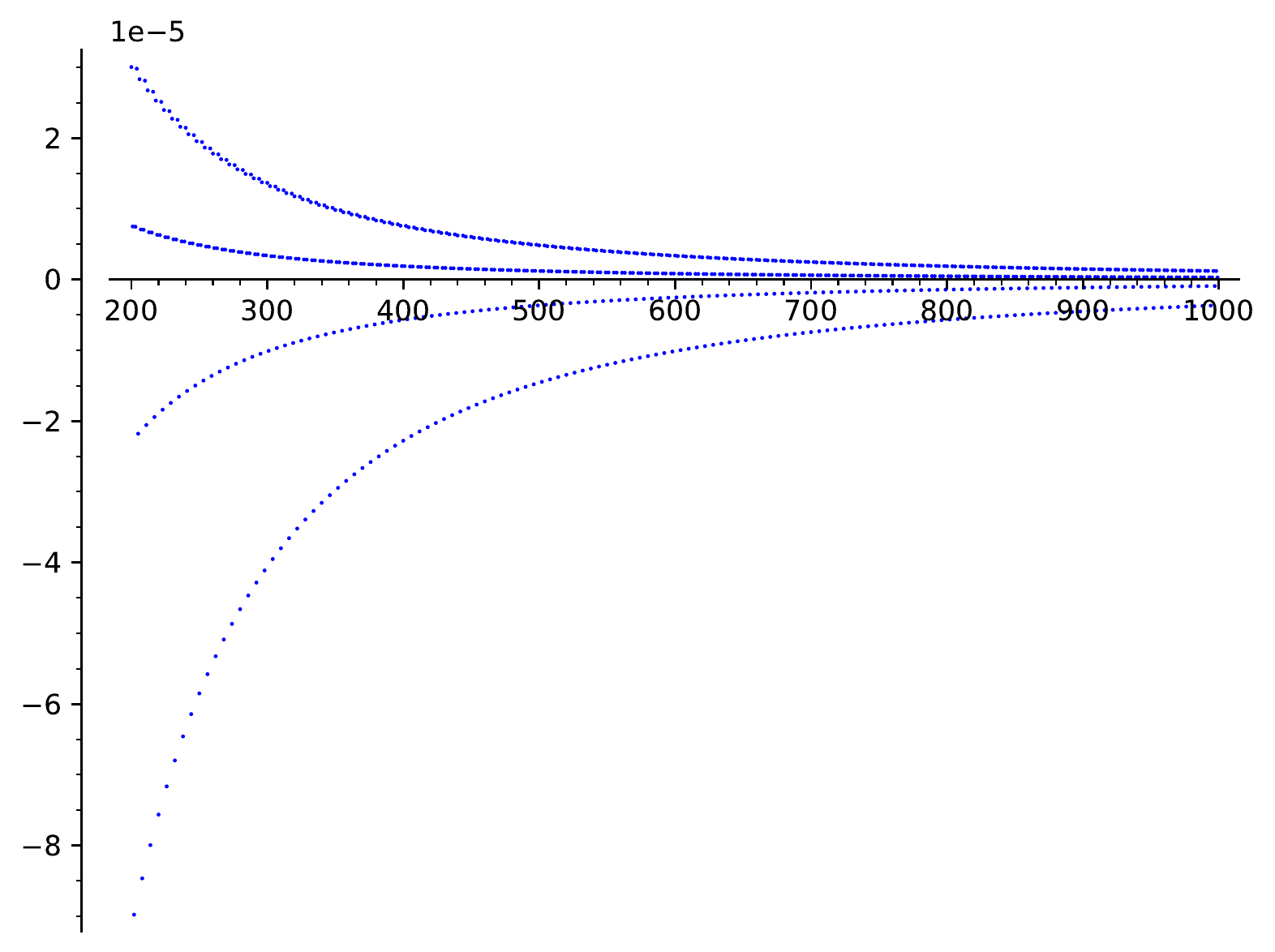}
		\caption{{\small $A_d = (1,d)$}}
	\end{subfigure} 
	\begin{subfigure}[t]{0.45\linewidth}
		\centering
		\includegraphics[width=\linewidth]{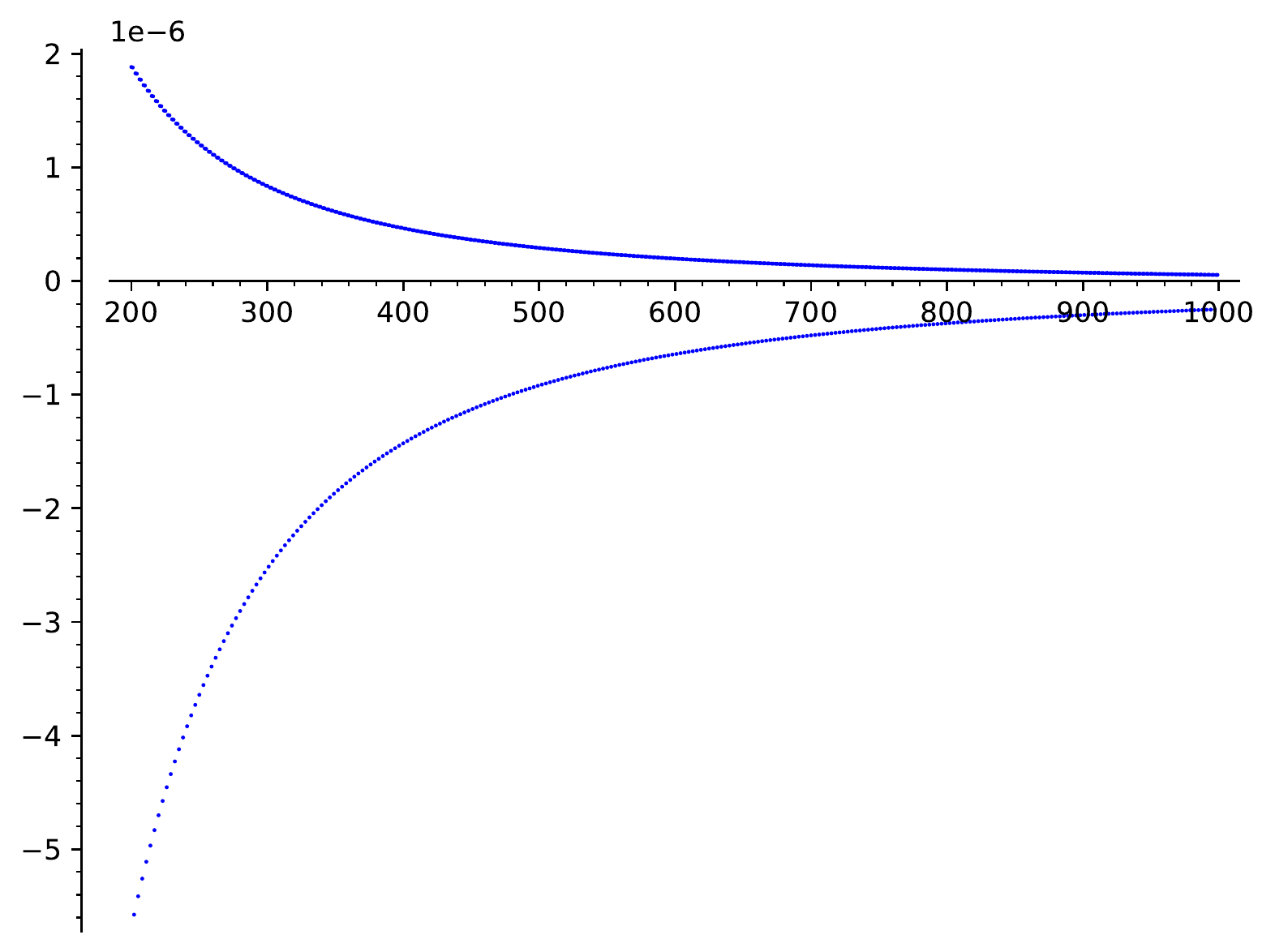}
		\caption{{\small $ A_d = (d,1) $}}
	\end{subfigure}
	
	\caption{Plots of $m(P_{A_d}) - m(P)$, for $P = z_1^2 + z_2 + 1$.}
	\label{fig:asymmetric}
\end{figure}

To conclude, we note that the invariant $\rho(A)$, whose divergence is sufficient to guarantee the convergence $m(P_A) \to m(P)$ (as we showed in \cref{thm:Lawton_generalisation}), will not suffice to express even the first term in the asymptotic expansion of $m(P_A) - m(P)$. More precisely, let $P = z_1^2 + z_2 + 1$, and consider the two sequences of matrices $A_d = (1,d)$ and $A_d = (d,1)$. Then $\rho(A_d) = d$ in both cases, but the convergence patterns for $m(P_{A_d}) \to m(P)$, portrayed in \cref{fig:asymmetric}, are quite different, which can be rigorously proved using Condon's formula \eqref{eq:Condon}. 

\section*{Acknowledgments}

The authors thank Vesselin Dimitrov and Philipp Habegger for fruitful discussions around their paper, and Wadim Zudilin for valuable comments on this work. The fourth named author would like to thank Chiara Amorino, Raphaël Ducatez, Roberto Gualdi and Tommaso Russo for useful discussions.

\section*{Funding}

The first and fourth named authors were supported by the research project ``Motivic homotopy, quadratic invariants and diagonal classes'' (ANR-21-CE40-0015) operated by the French National Research Agency (ANR).
This work was also performed within the framework of the LABEX MILYON (ANR-10-LABX-0070) of Université de Lyon, within the program "Investissements d'Avenir" (ANR-11-IDEX-0007), also operated by the French National Research Agency (ANR).

\section*{Addresses}
	
	François Brunault, \textsc{UMPA, École normale supérieure de Lyon, 46 allée d'Italie, 69100 Lyon, France}\par\nopagebreak
	\textit{E-mail address:} \href{mailto:francois.brunault@ens-lyon.fr}{\texttt{francois.brunault@ens-lyon.fr}}
	
	\medskip
	\noindent
	Antonin Guilloux, \textsc{Sorbonne Université, 4 place Jussieu, Boite Courrier 247, 75252 Paris Cedex 5, France} \par\nopagebreak
	\textit{E-mail address:} \href{mailto:antonin.guilloux@imj-prg.fr}{\texttt{antonin.guilloux@imj-prg.fr}}
	
	\medskip
	\noindent
	Mahya Mehrabdollahei, \textsc{Sorbonne Université, 4 place Jussieu, 75252 Paris Cedex 5, France} \par\nopagebreak
	\textit{E-mail address:} \href{mailto:mahya.mehrabdollahei@imj-prg.fr}{\texttt{mahya.mehrabdollahei@imj-prg.fr}}
	
	\medskip
	\noindent
	Riccardo Pengo, \textsc{UMPA, École normale supérieure de Lyon, 46 allée d'Italie, 69100 Lyon, France}\par\nopagebreak
	\textit{E-mail address:} \href{mailto:riccardo.pengo@ens-lyon.fr}{\texttt{riccardo.pengo@ens-lyon.fr}}
	

\begin{thebibliography}{9} 
		\bibitem{Akhtari_Vaaler_2019}
		Akhtari, S., \& Vaaler, J. D. (2019). \textit{Lower bounds for Mahler measure that depend on the number of monomials}. International Journal of Number Theory, 15(07), 1425–1436. doi:\href{https://doi.org/10.1142/S1793042119500805}{10.1142/S1793042119500805}. arXiv:\href{https://arxiv.org/abs/1810.12413}{1810.12413}.
		
		\bibitem{Apostol_1976}
		Apostol, T. M. (1976). \textit{Introduction to analytic number theory}. Springer-Verlag, New York-Heidelberg. doi:\href{https://doi.org/10.1007/978-1-4757-5579-4}{10.1007/978-1-4757-5579-4}
		
		\bibitem{Belabas_Cohen_2021}
		Belabas, K., \& Cohen, H. (2021). \href{https://www.math.u-bordeaux.fr/~kbelabas/Numerical_Algorithms/}{\textit{Numerical algorithms for number theory—using Pari/GP.}} Mathematical Surveys and Monographs, 254. American Mathematical Society, Providence, RI.
		
		\bibitem{Betke_Henk_Wills_1993}
		Betke, U., Henk, M., \& Wills, J. M. (1993). \textit{Successive-minima-type inequalities}. Discrete \& Computational Geometry, 9(2), 165–175. doi:\href{https://doi.org/10.1007/BF02189316}{10.1007/BF02189316}
		
		\bibitem{Bleistein_Handelsman_1986}
		Bleistein, N., \& Handelsman, R. A. (1986). \textit{Asymptotic expansions of integrals} (Second edition). Dover Publications, Inc., New York. ISBN:\href{https://store.doverpublications.com/0486650820.html}{978-0486650821}
		
		\bibitem{Boyd_1981}
		Boyd, D. W. (1981). \textit{Kronecker’s theorem and Lehmer’s problem for polynomials in several variables}. Journal of Number Theory, 13(1), 116–121. doi:\href{https://doi.org/10.1016/0022-314X(81)90033-0}{10.1016/0022-314X(81)90033-0} 
		
		\bibitem{Boyd_1981_Speculations}
		Boyd, D. W. (1981). \textit{Speculations concerning the range of Mahler’s measure}. Canadian Mathematical Bulletin, 24(4), 453–469. doi:\href{https://doi.org/10.4153/CMB-1981-069-5}{10.4153/CMB-1981-069-5}
		
		\bibitem{Boyd_1998}
		Boyd, D. W. (1998). \textit{Mahler’s measure and special values of $L$-functions.} Experimental Mathematics, 7(1), 37–82. doi:\href{https://doi.org/10.1080/10586458.1998.10504357}{10.1080/10586458.1998.10504357}
		
		\bibitem{Boyd_Mossinghoff_2005}
		Boyd, D. W., \& Mossinghoff, M. J. (2005). \textit{Small Limit Points of Mahler’s Measure}. Experimental Mathematics, 14(4), 403–414. doi:\href{https://doi.org/10.1080/10586458.2005.10128936}{10.1080/10586458.2005.10128936}
		
		\bibitem{Boyd_Rodriguez-Villegas_2001}
		Boyd, D. W., \& Rodriguez-Villegas, F. (2002). \textit{Mahler’s Measure and the Dilogarithm (I)}. Canadian Journal of Mathematics, 54(3), 468–492. doi:\href{https://doi.org/10.4153/CJM-2002-016-9}{10.4153/CJM-2002-016-9}
		
		\bibitem{Brunault_Zudilin_2020}
		Brunault, F., \& Zudilin, W. (2020). \textit{Many Variations of Mahler Measures: A Lasting Symphony.} Cambridge University Press. ISBN: \href{https://www.cambridge.org/core/books/many-variations-of-mahler-measures/29DB6CD1A87B356AD304DED9ECC9F4EE}{978-1-108-79445-9}
		
		\bibitem{Cai_Todd_2014}
		Cai, X., \& Todd, R. G. (2014). \textit{A cellular basis for the generalized Temperley–Lieb algebra and Mahler measure}. Topology and Its Applications, 178, 107–124. doi:\href{https://doi.org/10.1016/j.topol.2014.09.006}{10.1016/j.topol.2014.09.006}. arXiv:\href{https://arxiv.org/abs/1306.2585}{1306.2585}.
		
		\bibitem{Cartan_1961}
		Cartan, H. (1961). \textit{Théorie élémentaire des fonctions analytiques d’une ou plusieurs variables complexes}. Enseignement des Sciences. Hermann, Paris. ISBN:\href{https://www.editions-hermann.fr/livre/theorie-elementaire-des-fonctions-analytiques-d-une-ou-plusieurs-variables-complexes-henri-cartan}{978-2705652159}.
		
		\bibitem{Carter_Lalin_Manes_Miller_Mocz_2021}
		Carter, A., Lalín, M., Manes, M., Miller, A. B., \& Mocz, L. (2022). \textit{Two-variable polynomials with dynamical Mahler measure zero}. Research in Number Theory, 8(2), 25. doi:\href{https://doi.org/10.1007/s40993-022-00322-z}{10.1007/s40993-022-00322-z}. arXiv:\href{http://arxiv.org/abs/2110.06496}{2110.06496}. 
		
		
		\bibitem{Champanerkar_Kofman_2005}
		Champanerkar, A., \& Kofman, I. (2005). \textit{On the Mahler measure of Jones polynomials under twisting}. Algebraic \& Geometric Topology, 5(1), 1–22. doi:\href{https://doi.org/10.2140/agt.2005.5.1}{10.2140/agt.2005.5.1}. arXiv:\href{https://arxiv.org/abs/math/0404236}{math/0404236}.
		
		\bibitem{Champanerkar_Kofman_2006}
		Champanerkar, A., \& Kofman, I. (2006). \textit{On links with cyclotomic Jones polynomials}. Algebraic \& Geometric Topology, 6(4), 1655–1668. doi:\href{https://doi.org/10.2140/agt.2006.6.1655}{10.2140/agt.2006.6.1655}. arXiv:\href{https://arxiv.org/abs/math/0605631}{arXiv:math/0605631}.
		
		\bibitem{Condon_2012}
		Condon, J. D. (2012). \textit{Asymptotic expansion of the difference of two Mahler measures}. Journal of Number Theory, 132(9), 1962–1983. doi:\href{https://doi.org/10.1016/j.jnt.2012.02.022}{10.1016/j.jnt.2012.02.022}. arXiv:\href{https://arxiv.org/abs/1111.0241}{1111.0241}.

		\bibitem{D'Andrea_Lalin_2007}
		D’Andrea, C., \& Lalín, M. N. (2007). \textit{On the Mahler measure of resultants in small dimensions}. Journal of Pure and Applied Algebra, 209(2), 393–410. doi:\href{https://doi.org/10.1016/j.jpaa.2006.06.004}{10.1016/j.jpaa.2006.06.004}, arXiv:\href{https://arxiv.org/abs/math/0604359}{0604359}
		
		\bibitem{Deninger_2009}
		Deninger, C. (2009). \textit{Mahler measures and Fuglede-Kadison determinants}. Münster Journal of Mathematics, 2, 45–63. \\ url:\url{https://www.uni-muenster.de/FB10/mjm/vol_2/mjm_vol_2_04.pdf}
		
		\bibitem{Dimitrov_Habegger_2019}
		Dimitrov, V., \& Habegger, P. (2019). \textit{Galois orbits of torsion points near atoral sets}, arXiv:\href{https://arxiv.org/abs/1909.06051v2}{1909.06051v2}
		
		\bibitem{Dobrowolski_2012}
		Dobrowolski, E. (2012). \textit{On a question of Schinzel about the length and Mahler’s measure of polynomials that have a zero on the unit circle}. Acta Arithmetica, 155, 453–463. doi:\href{https://doi.org/10.4064/aa155-4-8}{10.4064/aa155-4-8}
		
		\bibitem{Dobrowolski_2017}
		Dobrowolski, E. (2017). \textit{A Note on Lawton’s Theorem}. Canadian Mathematical Bulletin, 60(3), 484–489. doi:\href{https://doi.org/10.4153/CMB-2016-066-x}{10.4153/CMB-2016-066-x}
		
		\bibitem{Dobrowolski_Smyth_2017}
		Dobrowolski, E., \& Smyth, C. (2017). \textit{Mahler measures of polynomials that are sums of a bounded number of monomials}. International Journal of Number Theory, 13(06), 1603–1610. doi:\href{https://doi.org/10.1142/S1793042117500907}{10.1142/S1793042117500907}. arXiv:\href{https://arxiv.org/abs/1606.04376}{1606.04376}.
		
		\bibitem{Dubickas_2018}
		Dubickas, A. (2018). \textit{On sums of two and three roots of unity}. Journal of Number Theory, 192, 65–79. doi:\href{https://doi.org/10.1016/j.jnt.2018.03.017}{10.1016/j.jnt.2018.03.017}
		
		\bibitem{Dubickas_Jankauskas_2013}
		Dubickas, A., \& Jankauskas, J. (2013). \textit{Nonreciprocal algebraic numbers of small Mahler’s measure}. Acta Arithmetica, 157, 357–364. doi:\href{https://doi.org/10.4064/aa157-4-3}{10.4064/aa157-4-3}
		
		\bibitem{Duke_2007}
		Duke, W. (2007). \textit{A combinatorial problem related to Mahler’s measure}. Bulletin of the London Mathematical Society, 39(5), 741–748. doi:\href{https://doi.org/10.1112/blms/bdm057}{10.1112/blms/bdm057}
		
		\bibitem{Erdelyi_1981}
		Erdélyi, A., Magnus, W., Oberhettinger, F., \& Tricomi, F. G. (1981). \textit{Higher transcendental functions. Vol. II}. Robert E. Krieger Publishing Co., Inc., Melbourne, Fla.
		
		\bibitem{Everest_Ward_1999}
		Everest, G., \& Ward, T. (1999). \textit{Heights of Polynomials and Entropy in Algebraic Dynamics}. Springer London. doi:\href{https://doi.org/10.1007/978-1-4471-3898-3}{10.1007/978-1-4471-3898-3}

		\bibitem{Gon_Oyanagi_2004}
		Gon, Y., \& Oyanagi, H. (2004). \textit{Generalized Mahler measures and multiple sine functions}. International Journal of Mathematics, 15(05), 425–442. doi:\href{https://doi.org/10.1142/S0129167X04002363}{10.1142/S0129167X04002363}

		\bibitem{Goncalves_1950}
		Goncalves, J. V. (1950). \textit{L’inégalité de W. Specht}. Universidade de Lisboa. Revista da Faculdade de Ciências. $2^a$ Série. A: Ciências Matemáticas, 1, 167–171.
		
		\bibitem{Gu_Lalin_2021}
		Gu, J., \& Lalín, M. (2021). \textit{The Mahler measure of a three-variable family and an application to the Boyd–Lawton formula}. Research in Number Theory, 7(1), 13. doi:\href{https://doi.org/10.1007/s40993-021-00237-1}{10.1007/s40993-021-00237-1}.
		
		\bibitem{Guilloux_2021}
		Guilloux, A. (2021). \textit{Experimentation around speed of convergence in Boyd-Lawton theorem for Mahler measure}, \url{https://gitlab.inria.fr/aguillou/mahler_condon}
		
		\bibitem{Guilloux_Marche_2021}
		Guilloux, A., \& Marché, J. (2021). \textit{Volume function and Mahler measure of exact polynomials}. Compositio Mathematica, 157(4), 809–834. doi : \href{https://doi.org/10.1112/S0010437X21007016}{10.1112/S0010437X21007016}. arXiv:\href{https://arxiv.org/abs/1804.01395}{1804.01395}.
		
		\bibitem{Habegger_2018}
		Habegger, P. (2018). \textit{The norm of Gaussian periods}, Q. J. Math. 69, No. 1, 153--182 doi:\href{https://doi.org/10.1093/qmath/hax028Zbl}{10.1093/qmath/hax028Zbl}. arXiv:\href{https://arxiv.org/abs/1611.07287}{1611.07287}.
		
		\bibitem{Hajli_2020}
		Hajli, M. (2020). \textit{A new formula for Mahler’s measure}. Functiones et Approximatio Commentarii Mathematici, 62(2), 165–170. doi:\href{https://doi.org/10.7169/facm/1753}{10.7169/facm/1753}. 
		
		\bibitem{Issa_Lalin_2013}
		Issa, Z., \& Lalín, M. (2013). \textit{A Generalization of a Theorem of Boyd and Lawton}. Canadian Mathematical Bulletin, 56(4), 759–768. doi:\href{https://doi.org/10.4153/CMB-2012-010-2}{10.4153/CMB-2012-010-2}. arXiv:\href{https://arxiv.org/abs/1203.5379}{1203.5379}.
		
		\bibitem{Kitano_Morifuji_Takasawa_2004}
		Kitano, T., Morifuji, T., \& Takasawa, M. (2004). \textit{$L^2$-torsion invariants of a surface bundle over $S^1$}. Journal of the Mathematical Society of Japan, 56(2), 503–518. doi:\href{https://doi.org/10.2969/jmsj/1191418642}{10.2969/jmsj/1191418642}
		
		\bibitem{Kurokawa_Lalin_Ochiai_2008}
		Kurokawa, N., Lalín, M., \& Ochiai, H. (2008). \textit{Higher Mahler measures and zeta functions}. Acta Arithmetica, 135, 269–297. doi:\href{https://doi.org/10.4064/aa135-3-5}{10.4064/aa135-3-5}. arXiv:\href{https://arxiv.org/abs/0908.0171}{0908.0171}.
		
		\bibitem{Lalin_Sinha_2011}
		Lalín, M., \& Sinha, K. (2011). \textit{Higher Mahler measure for cyclotomic polynomials and Lehmer’s question}. The Ramanujan Journal, 26(2), 257–294. doi:\href{https://doi.org/10.1007/s11139-010-9278-6}{10.1007/s11139-010-9278-6}. arXiv:\href{https://arxiv.org/abs/1106.1304}{1106.1304}.
		
		\bibitem{Lawton_1983}
		Lawton, W. M. (1983). \textit{A problem of Boyd concerning geometric means of polynomials.} Journal of Number Theory, 16(3), 356–362. doi:\href{https://doi.org/10.1016/0022-314X(83)90063-X}{10.1016/0022-314X(83)90063-X}
		
		\bibitem{Le_2014}
		Lê, T. T. Q. (2014). \textit{Homology torsion growth and Mahler measure}. Commentarii Mathematici Helvetici, 89(3), 719–757. doi:\href{https://doi.org/10.4171/cmh/332}{10.4171/cmh/332}. arXiv:\href{https://arxiv.org/abs/1010.4199}{1010.4199}.
		
		\bibitem{Lehmer_1933}
		Lehmer, D. H. (1933). \textit{Factorization of certain cyclotomic functions.} Annals of Mathematics. Second Series, 34(3), 461–479. doi:\href{https://doi.org/10.2307/1968172}{10.2307/1968172}
		
		\bibitem{Lind_Schmidt_Ward_1990}
		Lind, D., Schmidt, K., \& Ward, T. (1990). \textit{Mahler measure and entropy for commuting automorphisms of compact groups}. Inventiones Mathematicae, 101(1), 593–629. doi:\href{https://doi.org/10.1007/BF01231517}{10.1007/BF01231517}
		
		\bibitem{Liu_Qin_2021}
		Liu, H., \& Qin, H. (2021). \textit{Mahler Measure of Families of Polynomials Defining Genus $2$ and $3$ Curves.} Experimental Mathematics, 0(0), 1–16. doi:\href{https://doi.org/10.1080/10586458.2021.1926014}{10.1080/10586458.2021.1926014}. arXiv:\href{https://arxiv.org/abs/1910.10884}{1910.10884}.
		
		\bibitem{Luck_2015}
		Lück, W. (2015). \textit{Estimates for spectral density functions of matrices over $\mathbb{C}[\mathbb{Z}^d]$}. Annales Mathématiques Blaise Pascal, 22(1), 73–88. doi:\href{https://doi.org/10.5802/ambp.346}{10.5802/ambp.346}. arXiv:\href{https://arxiv.org/abs/1310.8564}{1310.8564}.
		
		\bibitem{Luck_2018}
		Lück, W. (2018). \textit{Twisting $L^2$-invariants with finite-dimensional representations}. Journal of Topology and Analysis, 10(04), 723–816. doi:\href{https://doi.org/10.1142/S1793525318500279}{10.1142/S1793525318500279}. arXiv:\href{https://arxiv.org/abs/1510.00057}{1510.00057}.
		
		\bibitem{Mahler_1962}
		Mahler, K. (1962). \textit{On Some Inequalities for Polynomials in Several Variables}. Journal of the London Mathematical Society, s1-37(1), 341–344. doi:\href{https://doi.org/10.1112/jlms/s1-37.1.341}{10.1112/jlms/s1-37.1.341}
		
		\bibitem{Mehrabdollahei_2021}
		Mehrabdollahei, M. (2021). \textit{Mahler measure of $P_d$ polynomials.} arXiv:\href{https://arxiv.org/abs/2101.07675v3}{2101.07675v3} 
		
		\bibitem{Mellit_2019}
		Mellit, A. (2019). \textit{Elliptic dilogarithms and parallel lines.} Journal of Number Theory, 204, 1–24. doi:\href{https://doi.org/10.1016/j.jnt.2019.03.019}{10.1016/j.jnt.2019.03.019}. arXiv:\href{https://arxiv.org/abs/1207.4722}{1207.4722}.
		
		\bibitem{Navot_1962}
		Navot, I. (1962). \textit{A Further Extension of the Euler-Maclaurin Summation Formula}. Journal of Mathematics and Physics, 41(1–4), 155–163. doi:\href{https://doi.org/10.1002/sapm1962411155}{10.1002/sapm1962411155}
		
		\bibitem{Otmani_Rhin_Sac-Epee_2019}
		Otmani, S. E., Rhin, G., \& Sac-Épée, J.-M. (2019). \textit{Finding New Limit Points of Mahler’s Measure by Genetic Algorithms}. Experimental Mathematics, 28(2), 129–131. doi:\href{https://doi.org/10.1080/10586458.2017.1357511}{10.1080/10586458.2017.1357511}
				
		\bibitem{Pierce_1916}
		Pierce, T. A. (1916). \textit{The numerical factors of the arithmetic forms $\prod_{i=1}^n (1\pm\alpha_i^m)$.} Annals of Mathematics. Second Series, 18(2), 53–64. doi:\href{https://doi.org/10.2307/2007169}{10.2307/2007169}

		\bibitem{Rudin}
		Rudin, W (1986). \textit{Real and Complex Analysis, Third edition.} McGraw-Hill Book Co., New York, 1987. xiv+416 pp.
		
		\bibitem{Raimbault_2012}
		Raimbault, J. (2012). \textit{Exponential growth of torsion in abelian coverings}. Algebraic \& Geometric Topology, 12(3), 1331–1372. doi:\href{https://doi.org/10.2140/agt.2012.12.1331}{10.2140/agt.2012.12.1331}. arXiv:\href{https://arxiv.org/abs/1012.3666}{1012.3666}.
		
		\bibitem{Schinzel_1997}
		Schinzel, A. (1997). \textit{On the Mahler measure of polynomials in many variables}. Acta Arithmetica, 79, 77–81. doi:\href{https://doi.org/10.4064/aa-79-1-77-81}{10.4064/aa-79-1-77-81}
		
		\bibitem{Silver_Williams_2004}
		Silver, D. S., \& Williams, S. G. (2004). \textit{Mahler Measure of Alexander Polynomials}. Journal of the London Mathematical Society, 69(3), 767–782. doi:\href{https://doi.org/10.1112/S0024610704005289}{10.1112/S0024610704005289}. arXiv:\href{https://arxiv.org/abs/math/0105234}{math/0105234}.
		
		\bibitem{Silver_Williams_2012}
		Silver, D. S., \& Williams, S. G. (2012). \textit{Twisted Alexander invariants of twisted links}. Journal of Knot Theory and Its Ramifications, 21(11), 1250118. doi:\href{https://doi.org/10.1142/S0218216512501180}{10.1142/S0218216512501180}. arXiv:\href{https://arxiv.org/abs/1202.1515}{1202.1515}.
		
		\bibitem{Smyth_1981}
		Smyth, C. J. (1981). \textit{On measures of polynomials in several variables.} Bulletin of the Australian Mathematical Society, 23(1), 49–63. doi:\href{https://doi.org/10.1017/S0004972700006894}{10.1017/S0004972700006894}
		
		\bibitem{Smyth_2002}
		Smyth, C. J. (2002). \textit{An explicit formula for the Mahler measure of a family of $3$-variable polynomials}. Journal de Théorie Des Nombres de Bordeaux, 14(2), 683–700. doi:\href{https://doi.org/10.5802/jtnb.382}{10.5802/jtnb.382}
		
		\bibitem{Smyth_2008}
		Smyth, C. J. (2008). \textit{The Mahler measure of algebraic numbers: A survey.} In: ``Number theory and polynomials'' (Vol. 352, pagg. 322–349). Cambridge Univ. Press, Cambridge. doi:\href{https://doi.org/10.1017/CBO9780511721274.021}{10.1017/CBO9780511721274.021}. arXiv:\href{https://arxiv.org/abs/math/0701397}{math/0701397}.
		
		\bibitem{Smyth_2018}
		Smyth, C. J. (2018). \textit{Closed sets of Mahler measures}. Proceedings of the American Mathematical Society, 146(6), 2359–2372. doi:\href{https://doi.org/10.1090/proc/13951}{10.1090/proc/13951}. arXiv:\href{https://arxiv.org/abs/1606.04338}{1606.04338}.
		
		\bibitem{Wiener_1932}
		Wiener, N. (1932). \textit{Tauberian Theorems}. Annals of Mathematics, 33(1), 1–100. doi:\href{https://doi.org/10.2307/1968102}{10.2307/1968102}
		
 	\end{thebibliography}
\end{document}